\documentclass[a4paper]{amsart}
\usepackage{amssymb}
\usepackage{mathtools}
\usepackage{mathrsfs}
\usepackage[hidelinks,bookmarksdepth=3]{hyperref}
\usepackage{tikz}
\usepackage{graphicx}
\usepackage{todonotes}
\usepackage[font=small]{caption,subcaption}
\usetikzlibrary{decorations.markings}
\allowdisplaybreaks

\tikzset{
	wing/.style={
		decoration={
			markings, mark=at position 0.8 with {\arrow{latex}}
		},
		postaction={decorate},
	},
}

\newcommand{\Z}{\mathbb{Z}}
\newcommand{\R}{\mathbb{R}}
\newcommand{\C}{\mathbb{C}}

\newcommand{\N}{\mathbb{N}}
\newcommand{\ii}{{\rm i}}
\newcommand{\cent}[1]{#1^\circ}

\newcommand{\graph}{\mathsf}

\newcommand{\gH}{\graph{H}}

\newcommand{\fF}{\mathfrak{F}}
\newcommand{\fG}{\mathfrak{G}}
\newcommand{\fR}{\mathfrak{R}}

\newcommand{\rotate}{\varsigma}

\newcommand{\re}{\operatorname{Re}}
\newcommand{\im}{\operatorname{Im}}
\newcommand{\Res}[2]{\operatorname{Res}\left(#1,#2\right)}

\newcommand{\conj}{\operatorname{conj}}

\theoremstyle{plain} %text of this environment is typesetted in italics
\newtheorem{theorem}{Theorem}
\newtheorem{proposition}{Proposition}
\newtheorem{lemma}[proposition]{Lemma}

\theoremstyle{definition} %text of this environment is typesetted in roman letters
\newtheorem{definition}{Definition}
\newtheorem{example}{Example}

\theoremstyle{remark}
\newtheorem{remark}{Remark}

\begin{document}

\title[Catenoid limits of saddle towers]{Catenoid limits of singly periodic minimal surfaces with Scherk-type ends}

\author{Hao Chen}
\address[Chen]{Institute of Mathematical Sciences, ShanghaiTech University, 201210 Pudong, Shanghai, China}
\email{chenhao5@shanghaitech.edu.cn}
\thanks{H.\ Chen was partially supported by Individual Research Grant from Deutsche Forschungsgemeinschaft within the project ``Defects in Triply Periodic Minimal Surfaces'', Projektnummer 398759432.}

\author{Peter Connor}
\address[Connor]{Department of Mathematical Sciences, Indiana University South Bend}
\email{pconnor@iusb.edu}

\author{Kevin Li}
\address[Li]{School of Science, Engineering, and Technology, Penn State Harrisburg}
\email{khl3@psu.edu}


\date{\today}

\maketitle
\begin{abstract}
	We construct families of embedded, singly periodic minimal surfaces of any
	genus $g$ in the quotient with any even number $2n>2$ of almost parallel
	Scherk ends.  A surface in such a family looks like $n$ parallel planes
	connected by $n-1+g$ small catenoid necks.  In the limit, the family
	converges to an $n$-sheeted vertical plane with $n-1+g$ singular points
	termed nodes in the quotient.  For the nodes to open up into catenoid necks,
	their locations must satisfy a set of balance equations whose solutions are
	given by the roots of Stieltjes polynomials.
\end{abstract}

The goal of this paper is to construct families of singly periodic minimal
surfaces (SPMSs) of any genus in the quotient with any even number $2n>2$ of
Scherk ends (asymptotic to vertical planes).  Each family is parameterized by a
small positive real number $\tau > 0$.  In the limit $\tau \to 0$, the Scherk
ends tend to be parallel, and the surface converges to an $n$-sheeted vertical
plane with singular points termed nodes.  As $\tau$ increases, the nodes open
up into catenoid necks, and the surface looks like parallel planes connected by
these catenoid necks.

There are many previously known examples of such SPMSs.  Scherk
\cite{scherk1835} discovered the examples with genus zero and four Scherk ends,
in 1835.  Karcher \cite{karcher1988} generalized Scherk's surface in 1988 with
any even number $2n>2$ of Scherk ends.  In this paper, examples of genus 0 will
be called ``Karcher--Scherk saddle towers'' or simply ``saddle towers'', and
saddle towers with four Scherk ends will be called ``Scherk saddle towers''.
Karcher also added handles between adjacent pairs of ends, producing SPMSs of
genus $n$ with $2n$ Scherk ends.  Traizet glued Scherk saddle towers into SPMSs
of genus $(n^2-3n+2)/2$ with $2n>2$ Scherk ends because he was desingularizing
simple arrangements of $n>1$ vertical planes.  In 2006, Mart\'{i}n and Ramos
Batista \cite{martin2006} replaced the ends of Costa surface by Scherk ends,
thereby constructed an embedded SPMS of genus one with six Scherk ends and, for
the first time, without any horizontal symmetry plane.  Hauswirth, Morabito,
and Rodr\'{i}guez \cite{hauswirth2009} generalized this result in 2009, using
an end-to-end gluing method to replace the ends of Costa--Hoffman--Meeks
surfaces by Scherk ends, thereby constructed SPMSs of higher genus with six
Scherk ends.  In 2010, da Silva and Ramos Batista \cite{dasilva2010}
constructed a SPMS of genus two with eight Scherk ends based on Costa surface.
Also, Hancco, Lobos, and Ramos Batista \cite{hancco2014} constructed SPMSs with
genus $2n$ and $2n$ Scherk ends in 2014.

The examples of da Silva and Ramos Batista as well as all examples of Traizet
admit catenoid limits that can be constructed using techniques in the present
paper.

\medskip

One motivation of this work is an ongoing project that addresses to various
technical details in the gluing constructions.

Roughly speaking, given any ``graph'' $G$ that embeds in the plane and
minimizes the length functional, one could desingularize $G \times \mathbb{R}$
into a SPMS by placing a saddle tower at each vertex.  Previously, this was
only proved for simple graphs under the assumption of an horizontal reflection
plane~\cite{traizet1996, traizet2001}.  Recently, we managed to allow the graph
to have parallel edges, to remove the horizontal reflection plane by Dehn
twist~\cite{saddle1}, and to prove embeddedness by analysing the bendings of
Scherk ends~\cite{saddle2}.

However, we still require that the vertices of $G$ are neither ``degenerate''
nor ``special''.  Here, a vertex of degree $2k$ is said to be \emph{degenerate}
(resp.\ \emph{special}) if $k$ (resp.\ $k-1$) of its adjacent edges extend in
the same direction while other $k$ (resp.\ $k-1$) edges extend in the opposite
direction.  This limitation is due to the fact that a saddle tower with $2k$
Scherk ends can not have $k-1$ ends extending in the same direction while other
$k-1$ ends extending in the opposite direction.  Therefore, it is not possible
to place a saddle tower at a degenerate or special vertex.

Nevertheless, we do know SPMSs that desingularize $G \times \mathbb{R}$ where
$G$ is a graph with degenerate vertex.  To include these in the gluing
construction, we need to place catenoid limits of saddle towers, as those
constructed in this paper, at degenerate vertices.  From this point of view,
the present paper can be seen as preparatory:  The insight gained here will
help us to glue saddle towers with catenoid limits of saddle towers in a future
project.

\medskip

This paper reproduces the main result of the thesis of the third named
author~\cite{li2012}.  Technically, the construction implemented
in~\cite{li2012} was in the spirit of~\cite{traizet2002}, which defines the
Gauss map and the Riemann surface at the same time, and the period of the
surface was assumed horizontal.  Here, for the convenience of future
applications, we present a construction in the spirit of~\cite{traizet2008,
saddle1, saddle2}, which defines all three Weierstrass integrands by
prescribing their periods, and the period of the surface is assumed vertical.
In particular, we will reveal that a balance condition in~\cite{li2012} is
actually a disguise of the balance of Scherk ends: The unit vectors in the
directions of the ends sum up to zero.

\subsection*{Acknowledgements}
We would like to thank the referee for many helpful comments.

\section{Main result}\label{sec:main}

\subsection{Configuration}

We consider $L+1$ vertical planes, $L\ge 1$, labeled by integers $l \in [1,
L+1]$.  Up to horizontal rotations, we assume that these planes are all
parallel to the $xz$-plane, which we identify to the complex plane $\C$, with
the $x$-axis (resp.\ $z$-axis) corresponding to the real (resp.\ imaginary) axis.
We use the term ``layer'' for the space between two adjacent parallel planes.
So there are $L$ layers.

We want $n_l \ge 1$ catenoid necks on layer $l$, i.e.\ between the planes $l$
and $l+1$, $1 \le l \le L$.  For convenience, we adopt the convention that $n_l
= 0$ if $l < 1$ or $l > L$, and write $N = \sum n_l$ for the total number of
necks.  Each neck is labeled by a pair $(l,k)$, where $1 \le l \le L$ and $1
\le k \le n_l$.

To each neck is associated a complex number $q_{l,k} \in \C^\times = \C
\setminus\{0\}$, $1 \le l \le L$, $1 \le k \le n_l$.  Then the positions of the
necks are prescribed at $\ln q_{l,k} + 2m\pi\ii$, $m \in \Z$.  Recall that the
$z$-axis is identified to the imaginary axis of the complex plane $\C$, so the
necks are periodic with period vector $(0,0,2\pi)$.  Note that, if we multiply
$q_{l,k}$'s by the same complex factor $c$, then the necks are all translated
by $\ln c \pmod{2\pi\ii}$.  So we may quotient out translations by fixing
$q_{1,1} = 1$. 

Moreover, each plane has two ends asymptotic to vertical planes.  We label the
end of plane $l$ that expands in the $-x$ (resp.\ $x$) direction by $0_l$
(resp.\ $\infty_l$).  To be compatible with the language of graph theory that
were used for gluing saddle towers~\cite{saddle1}, we use
\[
	\gH = \{\eta_l \colon 1 \le l \le L+1, \eta \in \{0,\infty\} \}
\]
to denote the set of ends.  When $0_l$ is used as subscript for
parameter $x$, we write $x_{l,0}$ instead of $x_{0_l}$ to ease the notation;
the same applies to $\infty_l$.

To each end is associated a real number $\dot\theta_h$, $h \in \gH$.  They
prescribe infinitesimal changes of the directions of the ends.  More
precisely, for small $\tau$, we want the unit vector in the direction of the
end $h$ to have a $y$-component or the order $\tau \dot\theta_h +
\mathcal{O}(\tau^2)$.  

\begin{remark}\label{rmk:normalize}
	Multiplying $\dot\theta$ by a common real constant leads to a
	reparameterization of the family.  Adding a common real constant to
	$\dot\theta_{l,0}$ and subtracting the same constant from
	$\dot\theta_{l,\infty}$ leads to horizontal rotations of the surface.
\end{remark}

In the following, we write
\[
	q = (q_{l,k})_{\substack{1 \leqslant l \leqslant L\\1 \leqslant k \leqslant n_l}}
	\quad\text{and}\quad
	\dot\theta = (\dot\theta_h)_{h \in \gH}.
\]
Then a \emph{configuration} refers to the pair $(q,\dot\theta)$.

% The catenoid necks on layer $l$ is given by $c_l$, recursively defined by
% \[
% 	c_l = \begin{cases}
% 		(\dot\theta_{1,0} + \dot\theta_{1,\infty})/n_1, & l = 1,\\
% 		(n_{l-1} c_{l-1}+\dot\theta_{l,0}+\dot\theta_{l,\infty})/n_l, & 2 \le l \le L.
% 	\end{cases}
% \]
% Note that
% \[
% 	- n_l c_l + n_{l-1} c_{l-1} + \dot\theta_{l,0} + \dot\theta_{l,\infty} = 0,\quad 1 \le l \le L+1.
% \]
% In particular, we have $c_L = -(\dot\theta_{L,0}+\dot\theta_{L,\infty}) / n_L$.

\subsection{Force}

Given a configuration $(q,\dot\theta)$, let $c_l$ be the real numbers that
solve 
\begin{equation}\label{eq:resc}
	- n_l c_l + n_{l-1} c_{l-1} + \dot\theta_{l,0} + \dot\theta_{l,\infty} = 0,
	\quad 1 \le l \le L+1.
\end{equation}
Recall the convention $n_l = 0$ if $l<1$ or $l>L$, so we also adopt the
convention $c_l = 0$ if $l<1$ or $l>L$. A summation over $l$ yields
\begin{equation}\label{eq:Theta1}
	\Theta_1 = \sum_{h \in \gH} \dot\theta_h = 0.
\end{equation}
If~\eqref{eq:Theta1} is satisfied, the real numbers $c_l$ are determined
by~\eqref{eq:resc} as functions of $\dot\theta$.

For $1 \le l \le L+1$, let $\psi_l$ be the meromorphic $1$-form on the Riemann
sphere $\hat\C$ with simple poles at $q_{l,k}$ with residue $-c_l$ for each $1
\le k \le n_l$, at $q_{l-1,k}$ with residue $c_{l-1}$ for each $1 \le k \le
n_{l-1}$, at $0$ with residue $\dot\theta_{l,0}$, and at $\infty$ with residue
$\dot\theta_{l,\infty}$.  More explicitly,
\[
	\psi_l =  \Big(\sum_{k = 1}^{n_l} \frac{-c_l}{z-q_{l,k}}
		+ \sum_{k = 1}^{n_{l-1}} \frac{c_{l-1}}{z-q_{l-1,k}} 
		+ \frac{\dot\theta_{l,0}}{z}
	\Big) dz.
\]
We then see that Equations~\eqref{eq:resc} arise from the Residue Theorem

\begin{remark}
  In the definition of configuration, we may replace $\dot\theta$ by the
  parameters $(c_l, \dot\theta_{l+1,0} - \dot\theta_{l,0})_{1 \le l \le L}$.
  Then $\dot\theta_{l,0}$'s are defined up to an additive constant
  (corresponding to a rotation), $\dot\theta_{l,\infty}$'s are determined
  by~\eqref{eq:resc}, and \eqref{eq:Theta1} is automatically satisfied.  To
  quotient out reparameterizations of the family, we may assume that $c_l = 1$
  for some $1 \le l \le L$.
\end{remark}

We define the force $F_{l,k}$ by
\begin{equation}\label{eq:force0}
	F_{l,k}
	% = \Res{q_{l,k}} \frac{\psi_l^2 + \psi_{l+1}^2}{d w_{l,k}}
	= \Res{\frac{\psi_l^2 + \psi_{l+1}^2}{2}\frac{z}{d z}}{q_{l,k}}.
\end{equation}
Or, more explicitly,
\begin{multline}\label{eq:force}
	F_{l,k} =
	\sum_{1 \le k \ne j \le n_l} \frac{2 c_l^2 q_{l,k}}{q_{l,k}-q_{l,j}}
	- \sum_{1 \le j \le n_{l+1}} \frac{c_l c_{l+1} q_{l,k}}{q_{l,k} - q_{l+1,j}} 
	- \sum_{1 \le j \le n_{l-1}} \frac{c_l c_{l-1} q_{l,k}}{q_{l,k} - q_{l-1,j}} \\
	+ c_l^2 + c_l (\dot\theta_{l+1,0} - \dot\theta_{l,0}).
\end{multline}
In~\cite{li2012}, the force had different formula depending on the parity of
$l$.  One verifies that both are equivalent to~\eqref{eq:force}.

\begin{remark}[Electrostatic interpretation]
	The force equation \eqref{eq:force} can be expressed as
	\begin{multline*}
		F_{l,k} =
		\sum_{1 \le k \ne j \le n_l} \frac{c_l^2 (q_{l,k}+q_{l,j})}{q_{l,k}-q_{l,j}}
		- \sum_{1 \le j \le n_{l+1}} \frac{c_l c_{l+1} (q_{l,k}+q_{l+1,j})}{2(q_{l,k} - q_{l+1,j})} \\
		- \sum_{1 \le j \le n_{l-1}} \frac{c_l c_{l-1} (q_{l,k}+q_{l-1,j})}{2(q_{l,k} - q_{l-1,j})} 
		+\frac{c_l}{2}(\dot\theta_{l,\infty}-\dot\theta_{l,0}-\dot\theta_{l+1,\infty}+\dot\theta_{l+1,0}).
	\end{multline*}

	%We have seen that both~\eqref{eq:kevinodd} and~\eqref{eq:kevineven} hold for
	%arbitrary $l$, adding them up yields
	%\begin{multline*}
		%2F_{l,k} =
		%\sum_{1 \le k \ne j \le n_l} 2 c_l^2 \frac{q_{l,k}+q_{l,j}}{q_{l,k}-q_{l,j}}\\
		%-\sum_{1 \le j \le n_{l+1}} c_l c_{l+1} \frac{q_{l,k} + q_{l+1,j}}{q_{l,k} - q_{l+1.j}}
		%-\sum_{1 \le j \le n_{l-1}} c_l c_{l-1} \frac{q_{l,k} + q_{l-1,j}}{q_{l,k} - q_{l-1.j}}\\
		%+c_l(\dot\theta_{l,\infty} - \dot\theta_{l,0} - \dot\theta_{l+1,\infty} + \dot\theta_{l+1,0}).
	%\end{multline*}
	Note that
	\[
		\frac{a+b}{a-b} = \coth\frac{\ln a - \ln b}{2}
		= \frac{2}{\ln a - \ln b} + \sum_{m=1}^{\infty}\Big(\frac{2}{\ln a - \ln b - 2m\pi\ii} + \frac{2}{\ln a - \ln b + 2m\pi\ii}\Big).
	\]
	Disregarding absolute convergence, we write this formally as
	\[
		\frac{a+b}{a-b} = \sum_{m \in \Z} \frac{2}{\ln a - \ln b - 2m\pi\ii}.
	\]
	Then the force is given, formally, by
	\begin{multline*}
		F_{l,k} =
		\sum_{0 \ne m \in \Z} \frac{2 c_l^2}{2m\pi\ii} 
		+\sum_{\substack{m \in \Z\\1 \le k \ne j \le n_l}} \frac{2 c_l^2}{\ln q_{l,k} - \ln q_{l,j} - 2m\pi\ii}\\
		-\sum_{\substack{m \in \Z\\1 \le j \le n_{l+1}}} \frac{c_l c_{l+1}}{\ln q_{l,k} - \ln q_{l+1.j} - 2m\pi\ii}
		-\sum_{\substack{m \in \Z\\1 \le j \le n_{l-1}}} \frac{c_l c_{l-1}}{\ln q_{l,k} - \ln q_{l-1.j} - 2m\pi\ii}\\
		+\frac{c_l}{2}(\dot\theta_{l,\infty}-\dot\theta_{l,0}-\dot\theta_{l+1,\infty}+\dot\theta_{l+1,0}).
	\end{multline*}
	Recall that $\ln q_{l,k} + 2m\pi\ii$ are the real positions of the necks.  So
	this formal expression has an electrostatic interpretation similar to those
	in~\cite{traizet2002} and~\cite{traizet2008}.  Here, each neck interacts not
	only with all other necks in the same or adjacent layers, but also with
	background constant fields given by $\dot\theta$.
\end{remark}

\begin{remark}[Another electrostatic interpretation]\label{rem:electrob}
	In fact, $\eqref{eq:force}/q_{l,k}$ has a similar electrostatic
	interpretation.  But this time, the necks are seen as placed at $q_{l,k}$.
	Each neck interacts with all other necks in the same and adjacent layers, as
	well as a virtual neck at $0$ with ``charge'' $c_l + \dot\theta_{l+1,0} -
	\dot\theta_{l,0}$.  This is no surprise, as electrostatic laws are known to be
	preserved under conformal mappings (such as $\ln(z)$).
\end{remark}

\subsection{Main result}

In the following, we write $F=(F_{l,k})_{1 \le l \le L, 1 \le k \le n_l}$.

\begin{definition}
	The configuration is \emph{balanced} if $F = 0$ and $\Theta_1=0$.
\end{definition}
Summing up all forces yields a necessary condition for the configuration to be
balanced, namely
\begin{multline*}
	\Theta_2 =\sum_{\substack{1 \le l \le L \\ 1 \le k \le n_l}} F_{l,k}
	= \sum_{1 \le l \le L+1} -\frac{1}{2} \Bigg(
		\Res{\frac{z\psi_l^2}{dz}}{0} + \Res{\frac{z\psi_l^2}{dz}}{\infty}
	\Bigg)\\
	= \sum_{1 \le l \le L+1} \frac{{\dot\theta_{l,\infty}^2 - \dot\theta_{l,0}^2}}{2} = 0.
\end{multline*}

\begin{lemma}\label{lem:rigid}
	The Jacobian matrix $\partial (\Theta_1, \Theta_2)/\partial \dot\theta$ has real
	rank $2$ as long as $c_l \ne 0$ for some $1 \le l \le L$.
\end{lemma}

The assumption of the lemma simply says that the surface does not remain a
degenerate plane to the first order.

\begin{proof}
	The proposition says that the matrix has an invertible minor of size $2
	\times 2$.  Explicitly, we have
	\[
		\frac{\partial (\Theta_1, \Theta_2)}{\partial(\dot\theta_{l,0}, \dot\theta_{l,\infty})}
		=\begin{pmatrix} 1 & 1\\ -\dot\theta_{l,0} & \dot\theta_{l,\infty} \end{pmatrix}.
	\]
	This minor is invertible if and only if the ends $\dot\theta_{l,0} +
	\dot\theta_{l,\infty} \ne 0$.  This must be the case for at least one $1 \le l
	\le L$ because, otherwise, we have $c_l = 0$ for all $1 \le l \le L$.
\end{proof}

\begin{definition}
	The configuration is \emph{rigid} if the complex rank of $\partial F
	/ \partial q$ is $N-1$.
\end{definition}

\begin{remark}
	In fact, the complex rank of $\partial F / \partial q$ is at most $N-1$.  We
	have seen that a complex scaling of $q$ corresponds to a translation of $\ln
	q_{l,k} + 2m\pi\ii$, $m \in \Z$, which does not change the force.  It then
	makes sense to normalize $q$ by fixing $q_{1,1}=1$.
\end{remark}

% {\color{blue}
% 	Computing from~\eqref{eq:force}, we may obtain two equivalent expressions for $W$, namely
% 	\[
% 		W = \sum_{l=1}^{L} n_l^2 c_l^2
% 		- \sum_{l=1}^{L-1} n_l n_{l+1} c_l c_{l+1}
% 		- \sum_{l=1}^{L-1} n_l c_l (\dot\theta_{l,0} - \dot\theta_{l+1,0})
% 		= -\sum_{l=1}^{L} n_l c_l (\dot\theta_{l,0} + \dot\theta_{l+1,\infty}),
% 	\]
% 	or, alternatively
% 	\[
% 		W 
% 		= \sum_{l=1}^{L} n_l^2 c_l^2
% 		- \sum_{l=2}^{L} n_l n_{l-1} c_l c_{l-1}
% 		- \sum_{l=1}^{L-1} n_l c_l (\dot\theta_{l,0} - \dot\theta_{l+1,0})
% 		= \sum_{l=1}^{L} n_l c_l (\dot\theta_{l,\infty} + \dot\theta_{l+1,0}).
% 	\]
% 	Taking average of the two yields
% 	\begin{align*}
% 		W &= \frac{1}{2}\sum_{k=1}^{L} n_lc_l(\dot\theta_{l,\infty} - \dot\theta_{l,0} +
% 		\dot\theta_{l+1,0} - \dot\theta_{l+1,\infty})
% 		= \frac{1}{2}\sum_{k=1}^{L+1} (n_lc_l-n_{l-1}c_{l-1})(\dot\theta_{l,\infty}-\dot\theta_{l,0})\\
% 		&= \frac{1}{2}\sum_{k=1}^{L+1} (\dot\theta_{l,\infty}+\dot\theta_{l,0})(\dot\theta_{l,\infty}-\dot\theta_{l,0}) = \sum_{k=1}^{L+1} \frac{\dot\theta_{l,\infty}^2 - \dot\theta_{l,0}^2}{2}
% 	\end{align*}
% }

\begin{theorem}
	Let $(q,\dot\theta)$ be a balanced and rigid configuration such that $c_l\neq 0$ for $1\leq l\leq L$.  Then for $\tau >
	0$ sufficiently small, there exists a smooth family $M_\tau$ of complete
	singly periodic minimal surface of genus $g=N-L$, period $(0,0,2\pi)$, and
	$2(L+1)$ Scherk ends such that, as $\tau \to 0$,
	\begin{itemize}
		\item $M_\tau$ converges to an $(L+1)$-sheeted $xz$-plane with singular
			points at $\ln q_{l,k} + 2m\pi\ii$, $m \in \Z$.  Here, the $xz$-plane is
			identified to the complex plane $\C$, with $x$-axis (resp.\ $z$-axis)
			identified to the real (resp.\ imaginary) axis.

		\item After suitable scaling and translation, each singular point opens up
			into a neck that converges to a catenoid.

		\item The unit vector in the direction of each Scherk end $h$ has the
			$y$-component $\tau\dot\theta_h + \mathcal{O}(\tau^2)$.
	\end{itemize}
	Moreover, $M_\tau$ is embedded if
	\begin{equation}\label{eq:embedcond}
		\dot\theta_{1,0} > \cdots > \dot\theta_{L+1,0}
		\quad\text{and}\quad
		\dot\theta_{1,\infty} > \cdots > \dot\theta_{L+1,\infty}.
	\end{equation}
\end{theorem}

\begin{remark}
	$M_\tau$ also depends smoothly on $\dot\theta$ belonging to the local smooth
	manifold defined by $\Theta_1=0$ and $\Theta_2 = 0$.  Up to
	reparameterizations of the family and horizontal rotations, we obtain
	families parameterized by $2L-1$ parameters.  Since we have $2(L+1)$ Scherk
	ends, this parameter count is compatible with the fact that Karcher--Scherk
	saddle towers with $2k$ ends form a family parameterized by $2k-3$
	parameters.
\end{remark}

\begin{remark}\label{rem:positivec}
	If the embeddedness condition~\eqref{eq:embedcond} is satisfied and
	$\Theta_1=0$, the sequence $\dot\theta_{l,0} + \dot\theta_{l,\infty}$ is
	strictly monotonically decreasing, and changes sign once and only once.  Then
	the sequence $n_l c_l$ is strictly concave (i.e.\ $n_{l-1} c_{l-1} + n_{l+1}
	c_{l+1} < 2 n_l c_l$ for $1 \le l \le L$).  Hence $c_l$, $1 \le l \le L$, are
	strictly positive, and the condition of Lemma~\ref{lem:rigid} is satisfied.
\end{remark}

\begin{remark}\label{rem:negativec}
We could allow some $c_l$ to be negative, with the price of losing embeddedness. Even worse, with negative $c_l$, the vertical planes in the limit will not be geometrically ordered as they are labeled.  For instance, if $L=2, c_1>0$, but $c_2<0$, then the catenoid necks, as well as the 1st and 3rd “planes", will all lie on the same side of the 2nd “plane".
\end{remark}

\begin{remark}\label{rem:zeroc}
We did not allow any $c_l$ to be 0 in Theorem 1.  Otherwise, the surface might still have nodes.  In that case, the claimed family might not be smooth, and the claimed genus would be incorrect.
\end{remark}
\section{Examples}

\subsection{Surfaces of genus 0}\label{ssec:g0}

When the genus $g = N-L = 0$, we have $n_l = 1$ for all $1 \le l \le L$, i.e.\
there is only one neck on every layer.  It then makes sense to drop the
subscript $k$.  For instance, the position and the force for the neck on layer
$l$ are simply denoted by $q_l$ and $F_l$.  We assume $L>1$ in this part.

In this case, if $\Theta_1 = 0$, Equation~\eqref{eq:resc} can be explicitly
solved by
\[
	c_l = \sum_{i=1}^{l} (\dot\theta_{i,0} + \dot\theta_{i,\infty}),\quad 1 \le l \le L,
\]
and the force can be written in the form 
\[
	F_{l} =
	- \widetilde{Q}_{l} + \widetilde{Q}_{l-1}
	+ c_l (\dot\theta_{l,\infty} + \dot\theta_{l+1,0}),\quad 1 \le l \le L,
\]
where we changed to the parameters
\[
	\widetilde{Q}_l = \frac{c_{l+1} c_l}{1-q_{l+1}/q_l}, \quad 1 \le l < L,
\]
with the convention that $\widetilde{Q}_0 = \widetilde{Q}_{L} = 0$.  Then the
forces are linear in $\widetilde Q$ and, if $\Theta_2 = 0$, the balance
condition $F=0$ is uniquely solved by
\begin{equation}\label{eq:summationQ}
	\widetilde{Q}_l
	= \sum_{i=1}^{l} c_i (\dot\theta_{i+1,0} + \dot\theta_{i,\infty})
	= - \sum_{i=l+1}^{L} c_i (\dot\theta_{i+1,0} + \dot\theta_{i,\infty}),
	\quad 1 \le l < L.
\end{equation}
Therefore, if we fix $q_1 = 1$, all other $q_l$, $1 < l \le L$ are uniquely
determined.

% Moreover, if the embeddedness condition~\eqref{eq:embedcond} is satisfied, the
% sequence $\dot\theta_{l,0} + \dot\theta_{l,\infty}$ is monotonically
% decreasing, and changes sign at most once.  Then the sequence $c$ is unimodal,
% i.e.\ there exists $1 \le l' \le L$ such that
% \[
% 	0 = c_0 < c_1 \le \cdots \le c_{l'} \ge \cdots \ge c_L > c_{L+1} = 0.
% \]

Recall from Remark~\ref{rem:positivec} that, under the embeddedness
condition~\eqref{eq:embedcond}, the numbers $c_l$, $1 \le l \le L$, are
positive.  Moreover, the summands in~\eqref{eq:summationQ} changes sign at most once, so the sequence $\widetilde Q$ is unimodal, i.e.\ there exists
$1 \le l' < L$ such that
\[
	0 = \widetilde{Q}_0 \le \widetilde{Q}_1 \le \cdots \le \widetilde{Q}_{l'} \ge \cdots \ge \widetilde{Q}_{L-1} \ge \widetilde{Q}_L = 0.
\]
Hence $\widetilde{Q}_l$, $1 \le l \le L$, are non-negative.  Moreover,
\begin{align*}
	\widetilde{Q}_l
	< \sum_{i=1}^{l} c_i (\dot\theta_{i,0} + \dot\theta_{i,\infty})
	&= \sum_{i=1}^{l} (c_i^2 - c_{i-1}c_i) \le c_l^2 \le c_{l+1}c_l
	& \text{if } l < l',\\
	\widetilde{Q}_l 
	< - \sum_{i=l+1}^{L} c_i (\dot\theta_{i+1,0} + \dot\theta_{i+1,\infty})
	&= \sum_{i=l+1}^{L} (c_i^2 - c_{i+1}c_i) \le c_{l+1}^2 \le c_{l+1}c_l
	& \text{if } l \ge l'.
\end{align*}
So $q$ consists of real numbers and $q_{l+1}/q_{l} < 0$ for all $1 \le l <
L$.

We have proved that
\begin{proposition}\label{prop:rigidhandel}
	If the genus $g = N-L = 0$, and $\dot\theta$ satisfies the balancing
	condition $\Theta_1 = \Theta_2 = 0$ as well as the embeddedness
	condition~\eqref{eq:embedcond}, then up to complex scalings, there exist
	unique values for the parameters $q$, depending analytically on $\dot\theta$,
	such that the configuration $(q,\dot\theta)$ is balanced.  Moreover, all such
	configurations are rigid.  If we fix $q_1 = 1$, then $q$ consist of real
	numbers, and we have $q_l > 0$ (resp.\ $<0$) if $l$ is odd (resp.\ even).
\end{proposition}

\subsection{Surfaces with four ends}

When $L=1$, $\Theta_1 = \Theta_2 = 0$ implies that
\[
	\dot\theta_{1,0} + \dot\theta_{2,\infty} =
	\dot\theta_{2,0} + \dot\theta_{1,\infty} = 0.
\]
Up to reparameterizations of the family, we may assume that $c_1 = 1$.  It
makes sense to drop the subscript $l$ in other notations, and write $F_k$ for
$F_{1,k}$, $q_k$ for $q_{1,k}$, and $n$ for $n_1$.  The goal of this part is to
prove the following classification result.

\begin{proposition}\label{prop:fourends}
	Up to a complex scaling, a configuration with $L=1$ and $n$ nodes must be
	given by $q_k = \exp(2\pi\ii k/n)$, and such a configuration is
	rigid.
\end{proposition}

Such a configuration is an $n$-covering of the configuration for Scherk saddle
towers.  As a consequence, the arising minimal surfaces are $n$-coverings of
Scherk saddle towers.  This is compatible with the result of~\cite{meeks2007}
that the Scherk saddle towers are the only connected SPMSs with four Scherk
ends.

\begin{proof}
	To find the positions $q_k$ such that
	\begin{equation}\label{eq:balance4ends}
		F_k=\sum_{1 \le k \ne j \le n} \frac{2q_k}{q_k-q_j} - (n-1) = 0,
		\qquad 1 \le k \le n,
	\end{equation}
	we use the polynomial method:  Consider the polynomial
	\[
		P(z) = \prod_{k = 1}^n (z - q_k).
	\]
	Then we have
	\begin{align*}
		P' &= P \sum_{k = 1}^n \frac{1}{z-q_k}\\
		P'' &= P \sum_{k = 1}^n \sum_{1 \le k \ne j \le n} \frac{1}{z-q_j}\frac{1}{z-q_k}
		= 2P \sum_{k = 1}^n \frac{1}{z-q_k} \sum_{1 \le k \ne j \le n} \frac{1}{q_k-q_j}\\
		&= P \sum_{k = 1}^n \frac{n-1}{q_k(z-q_k)}
		= (n-1) P \sum_{k=1}^n \frac{1}{z}\Big(\frac{1}{q_k} + \frac{1}{z-q_k}\Big)
		\qquad \text{by~\eqref{eq:balance4ends}}\\
		&= \frac{n-1}{z} \Big(P' - \frac{P'(0)}{P(0)} P \Big)
	\end{align*}
	For the last equation to have a polynomial solution, we must have $P'(0) =
	0$.  Otherwise, the left-hand side would be a polynomial of degree $n-2$, but
	the right-hand side would be a polynomial of degree $n-1$.

	Consequently, $F_k = 0$ if and only if
	\[
		zP''(z) -(n-1) P'(z) = 0,
	\]
	which, up to a complex scaling, is uniquely solved by
	\[
		P(z) = z^n - 1.
	\]
	So a balanced 4-end configuration must be given by the roots of unity $q_k =
	\exp(2\pi\ii k/n)$, $0 \le k \le n-1$.

	We now verify that the configuration is rigid.  For this purpose, we
	compute
	\[
		\frac{\partial F_k}{\partial q_j} = \begin{dcases}
			2 \frac{q_k}{(q_k-q_j)^2}, & j \ne k;\\
			2 \sum_{1 \le k \ne i \le n} \frac{-q_i}{(q_k-q_i)^2}, & j = k.
		\end{dcases}
	\]
	Note that $\sum_{j=1}^n q_j\partial F_k/\partial q_j = 0$ while
	\[
		q_j\frac{\partial F_k}{\partial q_j} = 2\frac{q_j q_k}{(q_k - q_j)^2}
		= 2\frac{e^{2\pi\ii \frac{j+k}{n}}}{(e^{2\pi\ii\frac{j}{n}} - e^{2\pi\ii\frac{k}{n}})^2}
		\in \R_{<0}
	\]
	when $j \ne k$, so the matrix
	\[
		\frac{\partial F}{\partial q} \operatorname{diag}(q_1, \cdots, q_n)
	\]
	has real entries, has a kernel of complex dimension $1$ (spanned by the
	all-one vector), and any of its principal submatrix is diagonally dominant.
	We then conclude that the matrix, as well as the Jacobian $\partial F /
	\partial q$, has a complex rank $n-1$.  This finishes the proof of
	rigidity.
\end{proof}

\begin{remark}
	The perturbation argument as in the proof of \cite[Proposition
	1]{traizet2002} also applies here, word by word, to prove the rigidity.
\end{remark}

\subsection{Gluing two saddle towers of different periods}

We want to construct a smooth family of configurations depending on a positive
real number $\lambda$ such that, for small $\lambda$, the configuration looks
like two columns of nodes far away from each other, one with period $2\pi/n_1$,
and the other with period $2\pi/n_2$.  If balanced and rigid, these
configurations would give rise to minimal surfaces that look like two Scherk
saddle towers with different periods that are glued along a pair of ends.  The
construction is in the same spirit as~\cite[\S~2.5]{traizet2002}
and~\cite[\S~4.3.4]{traizet2008}. 

\begin{proposition}
	For a real number $\lambda > 0$ sufficiently small, there are balanced and
	rigid configurations $(q(\lambda), \dot\theta(\lambda))$ with $L=2$
	depending smoothly on $\lambda$ such that, at $\lambda = 0$,
	\[
		\frac{q_{2,j}}{q_{1,k}} = 0, \qquad 1 \le k \le n_1, 1 \le j \le n_2.
	\]
	Up to a complex scaling and reparameterization, we may fix $q_{1,1}=1$, and
	write $q_{2,1}=\lambda\exp(\ii\phi)$.  Then, at $\lambda = 0$, we have
	\begin{equation}\label{eq:twoperiods}
		\dot\theta_{1,0} + \dot\theta_{2,\infty}
		=\dot\theta_{2,0} + \dot\theta_{3,\infty}
		=\dot\theta_{3,0} + \dot\theta_{1,\infty}
		= 0
	\end{equation}
	and, 
	\begin{align*}
		q_{1,k} &= \exp\Big(\frac{k-1}{n_1}2\pi\ii\Big),
		& 1 \le k \le n_1,\\
		\widetilde q_{2,k} := q_{2,k}/q_{2,1} &= \exp\Big(\frac{k-1}{n_2}2\pi\ii\Big),
		& 1 \le k \le n_2,
	\end{align*}
	where $\phi\operatorname{lcm}(n_1, n_2)$ is necessarily a
	multiple of $\pi$.
\end{proposition}

In other words, the construction only works if the configuration admits a
reflection symmetry.

\begin{remark}
	The first named author was shown a video suggesting that, when two Scherk
	saddle towers are glued into a minimal surface, one can slide one saddle
	tower with respect to the other while the surface remains minimal.  The
	proposition above suggests that this is not possible.
\end{remark}

In fact, the family of configurations also depends on $\dot\theta$ belonging to
the local manifold defined by $\Theta_1=\Theta_2=0$ and (one equation
from)~\eqref{eq:twoperiods}.  Up to rotations of the configuration and
reparameterizations of the family of minimal surfaces, the family of
configurations is parameterized, as expected, by two parameters.

\begin{proof}
	Let us first study the situation at $\lambda=0$.  We compute at $\lambda=0$
	\begin{align*}
		\frac{F_{1,k}}{c_1^2}
		=&\sum_{1\le k \ne j \le n_1} \frac{2q_{1,k}}{q_{1,k}-q_{1,j}}
		-\sum_{1 \le j \le n_2}\frac{c_2}{c_1}\frac{q_{1,k}}{q_{1,k}-q_{2,j}}
		+ 1 + \frac{\dot\theta_{2,0}-\dot\theta_{1,0}}{c_1}\\
		=& \sum_{1\le k \ne j \le n_1} \frac{2q_{1,k}}{q_{1,k}-q_{1,j}}
		- n_2 \frac{c_2}{c_1} + 1 + \frac{\dot\theta_{2,0}-\dot\theta_{1,0}}{c_1},\\
		\frac{F_{2,k}}{c_2^2}
		=&\sum_{1\le k \ne j \le n_2} \frac{2q_{2,k}}{q_{2,k}-q_{2,j}}
		-\sum_{1 \le j \le n_1} \frac{c_1}{c_2}\frac{q_{2,k}}{q_{2,k}-q_{1,j}}
		+ 1 + \frac{\dot\theta_{3,0}-\dot\theta_{2,0}}{c_2}\\
		=& \sum_{1\le k \ne j \le n_2} \frac{2q_{2,k}}{q_{2,k}-q_{2,j}}
		+ 1 + \frac{\dot\theta_{3,0}-\dot\theta_{2,0}}{c_2}.
	\end{align*}
	Write $G_l = \sum_k F_{l,k}$.  Summing the above over $k$ gives, at $\lambda=0$,
	\[
		\frac{1}{n_1}\frac{G_1}{c_1^2}=
		n_1 - n_2 \frac{c_2}{c_1} + \frac{\dot\theta_{2,0}-\dot\theta_{1,0}}{c_1},\qquad
		\frac{1}{n_2}\frac{G_2}{c_2^2} =
		n_2 + \frac{\dot\theta_{3,0}-\dot\theta_{2,0}}{c_2}.
	\]
	So $G_1=G_2=0$ at $\lambda=0$ only if 
	\begin{align*}
		0=&-(\dot\theta_{2,0}+\dot\theta_{3,\infty})=\dot\theta_{3,0}-\dot\theta_{2,0}+n_2c_2\\
		=&-(\dot\theta_{1,0}+\dot\theta_{2,\infty})=n_1 c_1 - n_2 c_2 + \dot\theta_{2,0}-\dot\theta_{1,0}.
	\end{align*}
	This together with $\Theta_1=0$ proves~\eqref{eq:twoperiods}.

	Now assume that~\eqref{eq:twoperiods} is satisfied.  Then we have, at $\lambda=0$
	\begin{align*}
		\frac{F_{1,k}}{c_1^2} =&
		\sum_{1 \le k \ne j \le n_1} \frac{2 q_{1,k}}{q_{1,k}-q_{1,j}} - (n_1 - 1),\\
		\frac{F_{2,k}}{c_2^2} =&
		\sum_{1 \le k \ne j \le n_2} \frac{2 q_{2,k}}{q_{2,k}-q_{2,j}} - (n_2 - 1).
	\end{align*}
	These expressions are identical to the force~\eqref{eq:balance4ends} for
	single layer configurations.  So we know for $l=1,2$ that, at $\lambda=0$,
	the configuration is balanced only if
	\[
		\widetilde{q}_{l,k}:=\frac{q_{l,k}}{q_{l,1}} = \exp\Big(\frac{k-1}{n_l}2\pi\ii\Big).
	\]
	Up to complex scaling, we may fix $q_{1,1}=1$ so $\widetilde{q}_{1,k} =
	q_{1,k}$.  And up to reparameterization of the family (of configurations), we
	write $q_{2,1}=\lambda\exp(\ii\phi)$.

	Now assume these initial values for $\widetilde{q}_{l,k}$.  Then we have, at
	$\lambda=0$,
	\begin{align*}
		\frac{G_2}{c_1 c_2} =& -\sum_{k=1}^{n_2}\sum_{j=1}^{n_1}\frac{q_{2,k}}{q_{2,k}-q_{1,j}}
		= \sum_{k=1}^{n_2}\sum_{j=1}^{n_1}\sum_{m=1}^{\infty} \Big(\frac{q_{2,k}}{q_{1,j}}\Big)^m\\
		=& \sum_{k=1}^{n_2}\sum_{j=1}^{n_1}\sum_{m=1}^{\infty}q_{2,1}^m \exp\bigg(2m\ii\pi\Big(\frac{k-1}{n_2} - \frac{j-1}{n_1}\Big)\bigg)
	\end{align*}
	Seen as a power series of $q_{2,1}$, the coefficient for $q_{2,1}^m$ is
	\[
		\sum_{k=1}^{n_2}\sum_{j=1}^{n_1}
		\exp\bigg(2m\ii\pi\Big(\frac{k-1}{n_2}-\frac{j-1}{n_1}\Big)\bigg).
	\]
	It is non-zero only if $m$ is a common multiple of $n_1$ and $n_2$, in which
	case the coefficient of $q_{2,1}^m$ equals $n_1n_2$. In particular, let $\mu
	= \operatorname{lcm}(n_1, n_2)$, then at $\lambda=0$,
	\begin{equation}\label{eq:imG2}
		\im \frac{G_2}{\lambda^{\mu}} = c_1 c_2 n_1 n_2 \sin(\mu\phi)
	\end{equation}
	vanishes if and only if $\mu\phi$ is a multiple of $\pi$.

	Now we use the Implicit Function Theorem to find balance configurations with
	$\lambda>0$.  From the proof for Proposition~\ref{prop:fourends}, we know
	that $\big(\frac{\partial F_{l,k}}{\partial \widetilde{q}_{l,j}}\big)_{2 \le
	j,k\le n_l}$, $l=1, 2$, are invertible.  Hence for $\lambda$ sufficiently
	small, there exist unique values for $(\widetilde{q}_{l,k})_{l = 1, 2; 2 \le
	k \le n_l}$, depending smoothly on $\lambda$, $\dot\theta$, and $\phi$, such
	that $(F_{l,k})_{l = 1, 2; 2 \le k \le n_l} = 0$.  By~\eqref{eq:imG2}, there
	exists a unique value for $\phi$, depending smoothly on $\lambda$ and
	$\dot\theta$, such that $\im G_2/\lambda^\mu = 0$.  Note also that $\re G_2$
	is linear in $\dot\theta$.  By Lemma~\ref{lem:rigid}, the solutions
	$(\lambda, \dot\theta)$ to $\re G_2 = 0$ and $\Theta_1 = \Theta_2 = 0$ form a
	manifold of dimension $4$ (including multiplication by common real factor on
	$\dot\theta$ and rotation of the configuration).  Finally, we have $G_1=0$ by
	the Residue Theorem, and the balance is proved.

	For the rigidity of the configurations with sufficiently small $\lambda$, we
	need to prove that the matrix
	\[
		\begin{pmatrix}
			\big(\frac{\partial F_{1,k}}{\partial q_{1,j}}\big)_{2 \le j,k\le n_1} & & \\
			& \big(\frac{\partial F_{2,k}}{\partial \widetilde{q}_{2,j}}\big)_{2 \le j,k\le n_2} & \\
			& & \frac{\partial G_2}{\partial q_{2,1}}
		\end{pmatrix}
	\]
	is invertible.  We know that the first two blocks are invertible at
	$\lambda=0$.  By continuity, they remain invertible for $\lambda$
	sufficiently small.  The last block is clearly non-zero for $\lambda \ne 0$
	sufficiently small.
\end{proof}

\subsection{Surfaces with six ends of type (n,1)} \label{ssec:n1}

In this section, we investigate examples with $L=2$ (hence six ends), $n_1 =
n$, $n_2 = 1$.  Up to a reparameterization of the family, we may assume that
$c_1 = 1$.  Up to a complex scaling, we may assume that $q_{2,1} = 1$.

We will prove that $q_{1,k}$'s are given by the roots of hypergeometric
polynomials.  Let us first recall their definitions.  A \emph{hypergeometric
function} is defined by 
\[
	{_2F_1}(a,b;c;z)=\sum_{k=0}^{\infty}\frac{(a)_k(b)_k}{(c)_k}\frac{z^k}{k!}
\]
with $a,b,c\in\C$, $c$ is not a non-positive integer,
\[
	(a)_k=a(a+1)\cdots(a+k-1)=\frac{\Gamma(a+k)}{\Gamma(a)},
\]
and $(a)_0 = 1$.  The hypergeometric function $w={_2F_1}(a,b;c;z)$ solves the
\emph{hypergeometric differential equation}
\begin{equation}\label{eq:ode2f1}
	z(1-z)w''+[c-(a+b+1)z]w'-abw=0.
\end{equation}

If $a=-n$ is a negative integer, 
\[
	{_2F_1}(-n, b; c; z) := \sum_{k = 0}^n (-1)^k \binom{n}{k} \frac{(b)_k}{(c)_k} z^k
\]
is a polynomial of degree $n$, and is referred to as a \emph{hypergeometric
polynomial}.

\begin{proposition} \label{thm:n1}
	Let $(q, \dot\theta)$ be a balanced configuration with $L=2$, $c_1 = 1$, $n_1
	= n$, $n_2 = 1$.  Then, up to a complex scaling, we have $q_{2,1} = 1$ and
	$(q_{1,k})_{1 \le k \le n}$ are the roots of the hypergeometric polynomial
	${_2F_1}(-n, b; c; z)$ with
	\[
		b := n - c_2 + \dot\theta_{2,0} - \dot\theta_{1,0},\qquad
		c := 1 + \dot\theta_{2,0} - \dot\theta_{1,0}.
	\]
	Moreover, as long as $b$ and $c$ are not non-positive integers, and $c-b$ is
	not a non-positive integer bigger than $-n$, the configuration is
	rigid.
\end{proposition}

\begin{proof}
	The force equations are
	\[
		\begin{split}
			F_{1,k}&=\sum_{1\le k \ne j \le n}\frac{2q_{1,k}}{q_{1,k}-q_{1,j}}-\frac{q_{1,k}c_2}{q_{1,k}-1}+c, 1\leq k\leq n,\\
			F_{2,1}&=-\sum_{j=1}^n\frac{c_2}{1-q_{1,j}}+c_2^2+c_2(\dot\theta_{3,0}-\dot\theta_{2,0}),
		\end{split}
	\]
	where $c := 1+\dot\theta_{2,0}-\dot\theta_{1,0}$.  To solve $F_{1,k}=0$ for
	$k=1,2,\ldots,n$, we use again the polynomial method: Let $P(z) =
	\prod_{k=1}^n (z-q_{1,k})$.  Then we have
	\begin{align*}
		P' &= P \sum_{k = 1}^n \frac{1}{z-q_{1,k}};\\
		\label{eq:ode6ends}
		P'' &= 2P \sum_{k = 1}^n \frac{1}{z-q_{1,k}} \sum_{1 \le k \ne j \le n} \frac{1}{q_{1,k}-q_{1,j}}\\
		&= P \sum_{k = 1}^n \frac{1}{(z-q_{1,k})}\Big(\frac{c_2}{q_{1,k}-1} - \frac{c}{q_{1,k}}\Big) \qquad \text{by $F_{1,k} = 0$}\\
		&= P \sum_{k=1}^n \Big( \frac{c_2}{(z-1)(z-q_{1,k})} + \frac{c_2}{(z-1)(q_{1,k}-1)}
			- \frac{c}{z (z-q_{1,k})}
		- \frac{c}{z q_{1,k}} \Big).
	\end{align*}
	So the configuration is balanced if and only if
	\begin{equation}\label{eq:ode6ends}
		P''
		+ \Big(\frac{-c_2}{z-1} + \frac{c}{z}\Big)P'
		+ \Big(\frac{c_2}{z-1}\frac{P'(1)}{P(1)}-\frac{c}{z}\frac{P'(0)}{P(0)}\Big)P=0.
	\end{equation}
	Define
	\[
		b := n - 1 - c_2 + c.
	\]
	For~\eqref{eq:ode6ends} to have a polynomial solution of degree $n$, we must
	have
	\[
		c_2\frac{P'(1)}{P(1)} = c \frac{P'(0)}{P(0)} = -nb,
	\]
	so that the leading coefficients cancel.  Then~\eqref{eq:ode6ends} becomes
	the hypergeometric differential equation
	\[
		z(1-z)P''+[c-(-n+b+1)z]P'+nbP=0
	\]
	to which the only polynomial solution (up to a multiplicative
	constant) is given by the hypergeometric polynomial $P(z) = {_2F_1}(-n, b; c;
	z)$ of degree $n$.

	Moreover, in order for $F_{2,1}=0$, we must have
	\begin{equation} \label{eq:F21}
		\dot\theta_{3,0}-\dot\theta_{2,0}=\sum_{j=1}^n\frac{1}{1-q_{1,j}}-c_2%=\frac{n}{2}.
		= \frac{P'(1)}{P(1)} - c_2 = -\frac{nb}{c_2} - c_2
	\end{equation}

	Note that $b$ and $c$ are real.  If $b$ is not a non-positive integer, and
	$c-b$ is not a non-positive integer bigger than $-n$, then all the $n$ roots
	of $P(z) = {_2F_1}(-n,b;c;z)$ are simple.  Indeed, under these assumptions,
	we have $P(0) = 1$ and $P(1) = (c-b)_n/(c)_n \ne 0$ by the Chu--Vandermonde
	identity.  Let $z_0$ be a root of $P(z)$, then $z_0 \ne 0, 1$.  In view of
	the hypergeometric differential equation, if $z_0$ is not simple, we have
	$P(z_0) = P'(z_0) = 0$ hence $P(z) \equiv 0$ by uniqueness theorem.

	The rigidity means that no perturbation of $q_{1,k}$ preserve the
	balance to the first order.  To prove this fact, we use a perturbation
	argument similar to that in the proof of \cite[Proposition 1]{traizet2002}.

	Let $(q_{1,k}(t))_{1 \le k \le n}$ be a deformation of the configuration such
	that $q_{1,k}(0) = q_{1,k}$ and $(\dot F_{1,k}(0))_{1 \le k \le n} = 0$,
	where dot denotes derivative with respect to $t$.  Define
	\[
		P_t(z) = \sum_{j = 0}^n a_j(t) z^j := \prod_{k = 1}^n (z-q_{1,k}(t)).
	\]
	Then we have
	\[
		z(1-z)P_t''+[c-(-n+b+1)z]P_t'+nbP_t = o(t),
	\]
	meaning that the coefficients from the left side are all $o(t)$.  So the
	coefficients of $P_t$ must satisfy
	\begin{equation}\label{eq:hypergeometriccoeff}
		% (-j(j-1)-(-n+b+1) j + nb) a_j(t) + (j(j+1) + cj) a_{j+1} = o(t).
		% (-j^2 + nj - bj + nb) a_j(t) + (j^2 + j + cj) a_{j+1}(t) &= o(t), & 0 \le j < n;\\
		(b+j)(n-j) a_j(t) + (j^2 + j + cj) a_{j+1}(t) = o(t), \qquad 0 \le j \le n.
	\end{equation}
	Note that $P_t(z)$ is monic by definition, meaning that $a_n(t) \equiv 1$.
	Since $b$ and $c$ are not non-positive integers, we conclude that $a_j(t) =
	o(t)$ for all $0 \le j \le n$.  The simple roots depend analytically on the
	coefficients, so $q_{1,k}(t) = q_{1,k} + o(t)$.
\end{proof}

The simple roots of ${_2F_1}(-n,b;c;z)$ are either real or form conjugate pairs.  As
a consequence, if rigid, the configurations in the proposition above
will give rise to minimal surfaces with horizontal symmetry planes.

\begin{example}
	For each integer $n \ge 2$, the real parameters $(b, c)$ for which
	${_2F_1}(-n, b; c; z)$ has only real simple roots has been enumerated
	in~\cite{dominici2013}.  The results are plotted in blue in
	Figure~\ref{fig:bcplot}.  The embeddedness conditions \eqref{eq:embedcond}
	are
	\begin{align*}
		\dot\theta_{1,0} > \dot\theta_{2,0} & \Rightarrow c < 1, \\
		\dot\theta_{1,\infty} > \dot\theta_{2,\infty} & \Rightarrow b > -n, \\
		\dot\theta_{2,0} > \dot\theta_{3,0} & \Rightarrow c_2^2 > -nb, \\
		\dot\theta_{2,\infty} > \dot\theta_{3,\infty} & \Rightarrow c_2^2 > n(c_2+b),
	\end{align*}
	where $c_2 = n-1-b+c$.  The region defined by these is plotted in red in
	Figure~\ref{fig:bcplot}.  Then non-integer parameters $(b,c)$ in the
	intersection of red and blue regions give rise to balanced and rigid
	configurations with real $q_{1,k}$.

\begin{figure}
	\includegraphics[width=\textwidth]{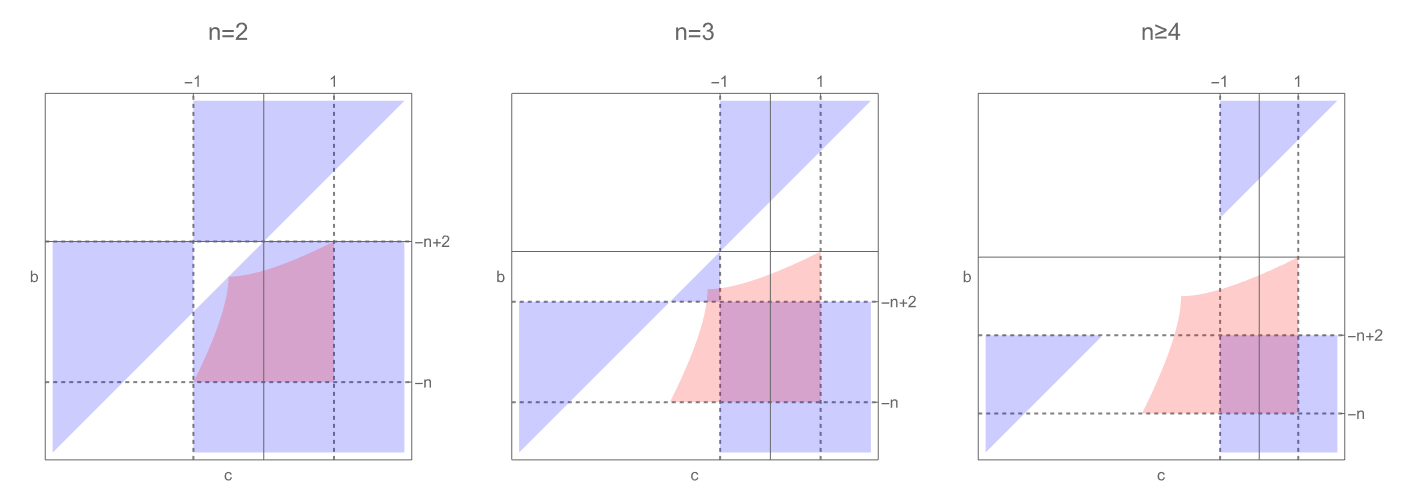}
	\caption{
		The set of $b$ and $c$ for which ${_2F_1}(-n,b;c;z)$ has only real simple
		roots (blue), and for which the embeddedness conditions are satisfied
		(red).\label{fig:bcplot}
	}
\end{figure}

	Figure \ref{fig:51a} shows the configurations of three examples with $n=5$. \qed
	
\begin{figure}[h]
  \centering
  \begin{subfigure}[b]{0.48\textwidth}
    \centering
    \includegraphics[width=\textwidth, trim = 0cm -.05cm 0cm 0cm]{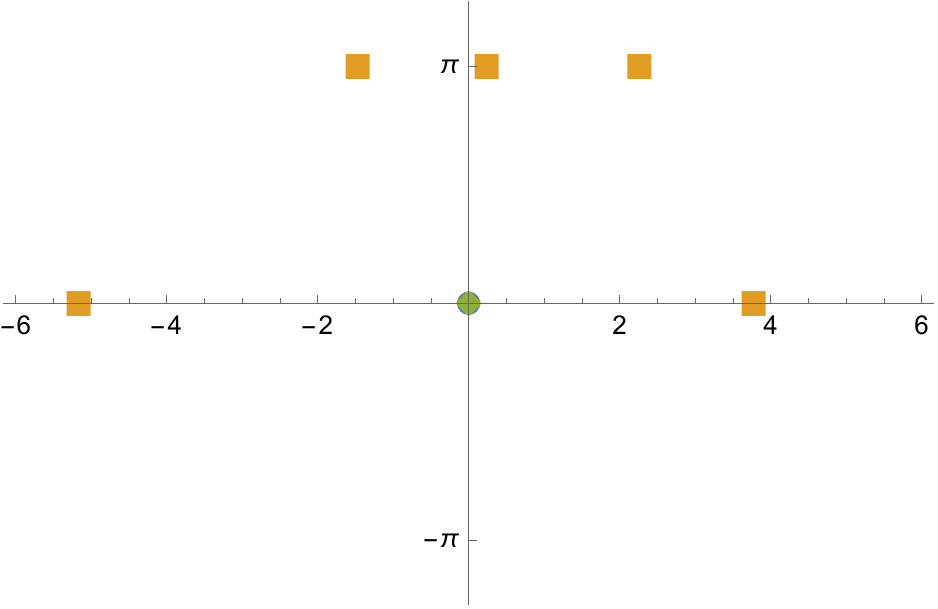}
    \caption*{
    	$b=-3.4, c=-0.1$}
  \end{subfigure} \hspace{.1in}
  ~ %add desired spacing between images, e. g. ~, \quad, \qquad, \hfill etc. 
  \begin{subfigure}[b]{0.48\textwidth}
    \centering
    \includegraphics[width=\textwidth]{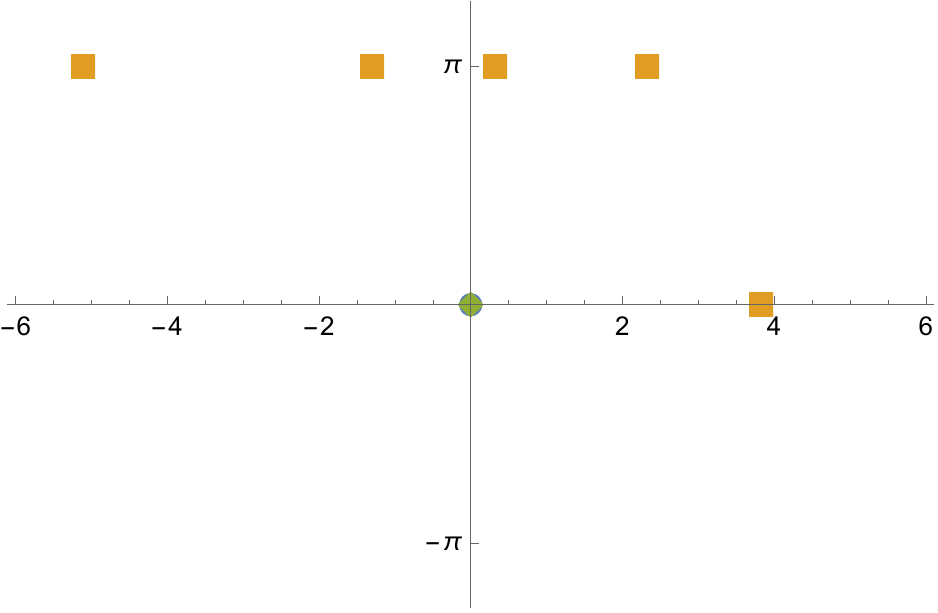}
    \caption*{
    	$b=-3.4, c=0.1$}
  \end{subfigure}\\
  
  ~ %add desired spacing between images, e. g. ~, \quad, \qquad, \hfill etc. 
  \begin{subfigure}[b]{0.65\textwidth}
    \centering
    \includegraphics[width=\textwidth]{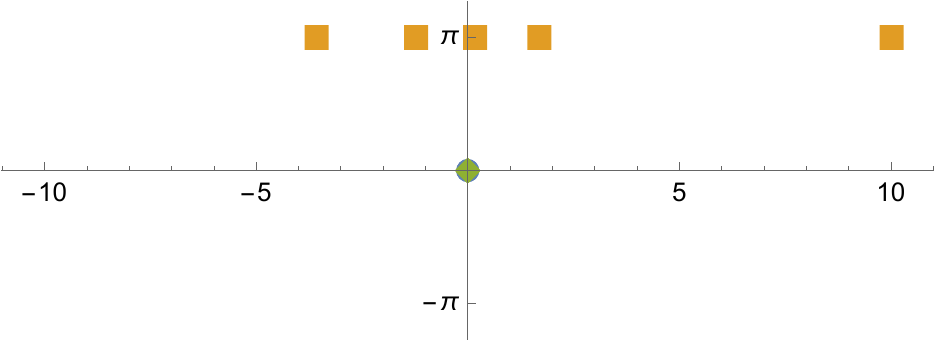}
    \caption*{
    	$b=-4.001, c=0.5$}
  \end{subfigure}
  ~ %add desired spacing between images, e. g. ~, \quad, \qquad, \hfill etc. 
 % \begin{subfigure}[b]{0.5\textwidth}
 %   \centering
 %   \includegraphics[width=\textwidth]{51b=-45c999.pdf}
 %   \caption{
%    	$b=-4.5, c=0.999$}
%  \end{subfigure}
  \caption{
  	$(5,1)$ balanced configurations.  The circles and squares represent the necks at levels one and two respectively.
	}\label{fig:51a}
\end{figure}
\end{example}

\begin{remark}
	As $c \to 0$, ${_2F_1}(-n,b;c;z)/\Gamma(c)$ converges to a polynomial with a
	root at $0$.  One may interpret that, as $c$ increases across $0$, a root
	moves from the interval $(-\infty, 0)$ to the interval $(0,1)$ through $0$.

	When $b = 1-n$, ${_2F_1}(-n,b;c;z)$ becomes a polynomial of degree $n-1$.
	One may interprete that, as $b$ increases across $1-n$, a root moves from the
	interval $(-\infty,0)$ to the interval $(1,\infty)$ through the infinity.
\end{remark}

\begin{example}\label{ex:2}
	Assume that $b+c=1-n$ (hence $c_2 = -2b$). Then by the identity
	\[
		{_2F_1}(-n,b;c;z)=\frac{(b)_n}{(c)_n}(-z)^n{_2F_1}\left(-n,1-c-n;1-b-n;\frac{1}{z}\right),
	\]
	the simple roots must be symmetrically placed.  That is, if $z_0$ is a
	root, so is $1/z_0$.  This symmetry appears in the resulting minimal surfaces
	as a rotational symmetry.  If the simple roots are real, the rotation reduces
	to a vertical reflectional.  In view of Figure~\ref{fig:bcplot}, we obtain
	the following concrete examples.
	\begin{itemize}
		\item $n \ge 2$ and $0 < c < 1$. In this case ${_2F_1}(-n, b; c; z)$ has
			$n$ simple negative roots.  See Figure \ref{fig:21} for an example of
			this type	with $n=2$.  Figure \ref{fig:51b} shows the configurations of two examples with $n=5$.

		\item $n \ge 3$ and $-1 < c < 0$, or $n=3$ and $-5/4<c<-1$, or $n
			= 2$ and $-1/2 < c < 0$. In these case, ${_2F_1}(-n, b; c; z)$ has $n-2$
			simple negative roots, one root $0 < z_0 < 1$, and another root $1/z_0 >
			1$.  	Figure \ref{fig:51c} shows the configurations of two examples with $n=5$. \qed
	\end{itemize}

	\begin{figure}[htbp] %%%%%%%%%%
		\centering
 		\includegraphics[height=1.7in]{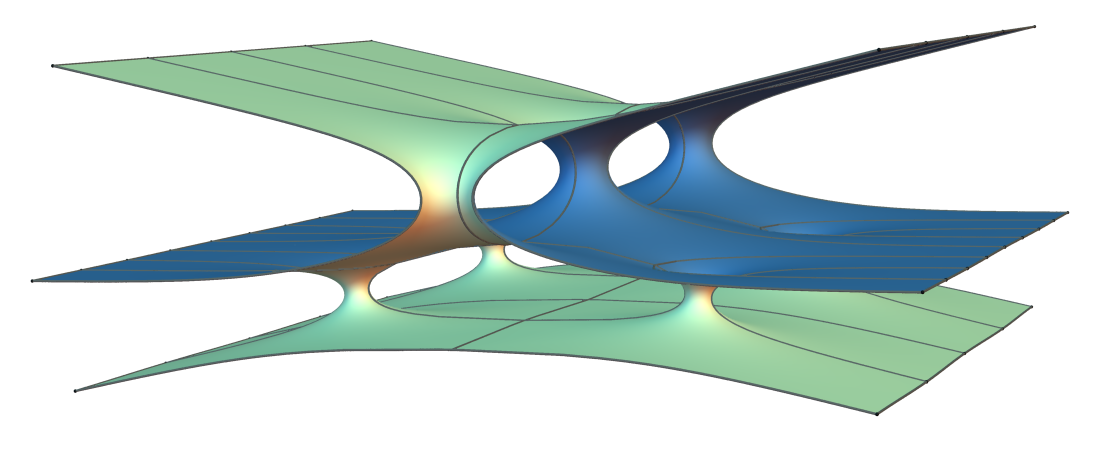}
		\caption{Genus 1 example with $n=2$ and $0<c<1$.}
		\label{fig:21}
	\end{figure} %%%%%%%%%%

\begin{figure}[h]
  \centering
  \begin{subfigure}[b]{0.48\textwidth}
    \centering
    \includegraphics[width=\textwidth, trim = 0cm -3.8cm 0cm 0cm]{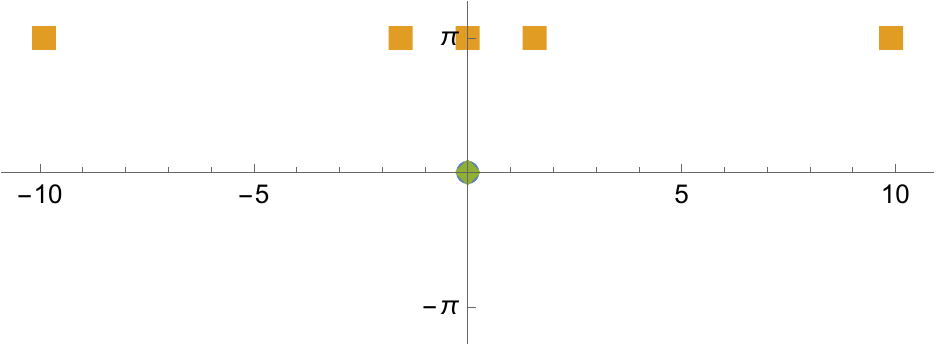}
    \caption*{
    	$c=0.001$}
  \end{subfigure} \hspace{.1in}
  ~ %add desired spacing between images, e. g. ~, \quad, \qquad, \hfill etc. 
  \begin{subfigure}[b]{0.48\textwidth}
    \centering
    \includegraphics[width=\textwidth]{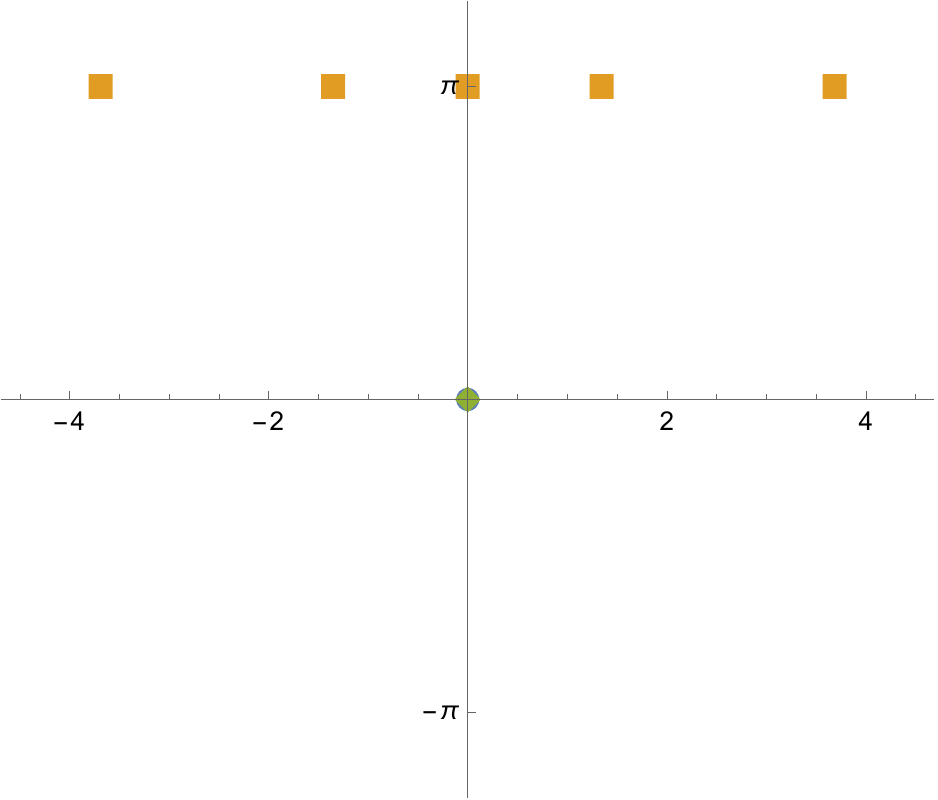}
    \caption*{
    	$c=0.5$}
  \end{subfigure}\\
  
  ~ %add desired spacing between images, e. g. ~, \quad, \qquad, \hfill etc. 
  \caption{
  	$(5,1)$ balanced configurations with $0<c<1$.  The circles and squares represent the necks at levels one and two respectively.
	}\label{fig:51b}
\end{figure}

\begin{figure}[h]
  \centering
  \begin{subfigure}[b]{0.48\textwidth}
    \centering
    \includegraphics[width=\textwidth]{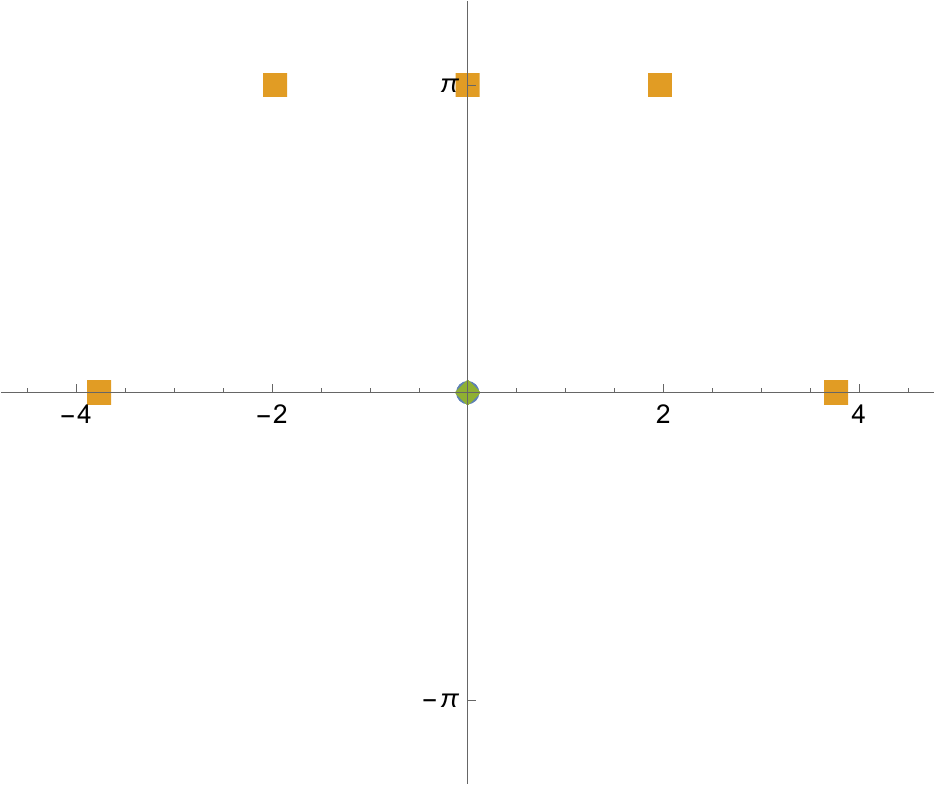}
    \caption*{
    	$c=-0.5$}
  \end{subfigure} \hspace{.1in}
  ~ %add desired spacing between images, e. g. ~, \quad, \qquad, \hfill etc. 
  \begin{subfigure}[b]{0.48\textwidth}
    \centering
    \includegraphics[width=\textwidth, trim = 0cm -1.8cm 0cm 0cm]{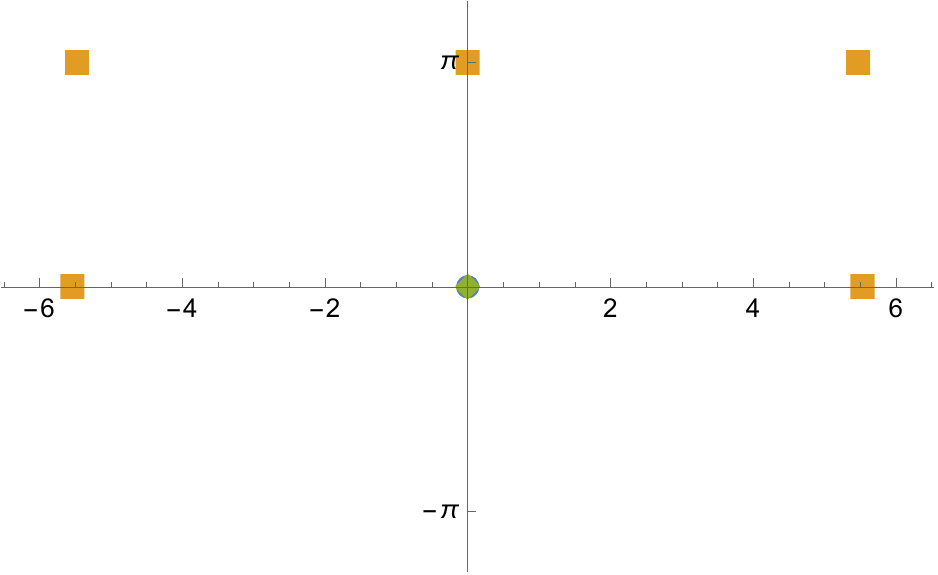}
    \caption*{
    	$c=-0.999$}
  \end{subfigure}\\
  
  ~ %add desired spacing between images, e. g. ~, \quad, \qquad, \hfill etc. 
  \caption{
  	$(5,1)$ balanced configurations with $-1<c<0$.  The circles and squares represent the necks at levels one and two respectively.
	}\label{fig:51c}
\end{figure}

\end{example}

\begin{remark}
	Examples with six Scherk ends are parameterized by three real parameters,
	here by $b$, $c$, and the family parameter $\tau$.  We see that the relation
	$b+c=1-n$ imposes a rotational symmetry.  It can be imagined that removing
	the relation would break this symmetry.
\end{remark}

\begin{remark}
	The polynomial method is often used to find balanced configurations of
	interacting points in the plane.  In minimal surface theory, it has been
	employed in many implementations of Traizet's node-opening technique
	\cite{traizet2002, traizet2002b, traizet2005, li2012, connor2012, connor2017,
	connor2017b, chen2022}.
	% Examples of singly periodic minimal surfaces are constructed in
	% \cite{connor2017b} with two annular ends and $n+1$ catenoid ends where the
	% location of the catenoid ends are the roots of a hypergeometric polynomial.
	% Examples of doubly periodic minimal surfaces are constructed in
	% \cite{connor2017}, in which balanced configurations are given by Heun
	% polynomials.
\end{remark}

\subsection{Surfaces with eight ends of type (1,n,1)} \label{ssec:1n1}
Theorem~\ref{thm:n1} generalizes to the following lemma with similar proof

\begin{lemma} \label{lemma:1n1}
	If we fix $q_{l\pm 1,k}$'s and assume that $c_l = 1$.  Then $q_{l,k}$'s in a
	balanced configuration are given by the roots of a Stieltjes polynomial
	$P(z)$ of degree $n_l$ that solves the generalized Lam\'e equation (a.k.a.\
	second-order Fuchsian equation)~\cite{marden1966}
	\begin{align}
		P''
		&+ \Big(
			\frac{c}{z}
			+ \sum_{k = 1}^{n_{l-1}} \frac{-c_{l-1}}{z-q_{l-1,k}}
			+ \sum_{k = 1}^{n_{l+1}} \frac{-c_{l+1}}{z-q_{l+1,k}}
		\Big)P'\label{eq:lame}\\
		&+ \Big(
			\frac{\gamma_0}{z}
			+ \sum_{k=1}^{n_{l-1}} \frac{\gamma_{l-1,k}}{z-q_{l-1,k}}
			+ \sum_{k=1}^{n_{l+1}} \frac{\gamma_{l+1,k}}{z-q_{l+1,k}}
		\Big)P =0,\nonumber
	\end{align}
	where $c = 1 + \dot\theta_{l+1,0} - \dot\theta_{l,0}$,
	% \[
	% 	\gamma_0 = -c \frac{P'(0)}{P(0)}
	% 	\quad\text{and}\quad
	% 	\gamma_{l\pm 1, k} = c_{l\pm 1} \frac{P'(q_{l\pm 1,k})}{P(q_{l\pm 1,k})},
	% \]
	subject to conditions
	\begin{align*}
		\gamma_0
		+ \sum_{k=1}^{n_{l-1}} \gamma_{l-1,k}
		+ \sum_{k=1}^{n_{l+1}} \gamma_{l+1,k}
		&=0,\\
		\sum_{k=1}^{n_{l-1}} \gamma_{l-1,k} q_{l-1,k}
		+ \sum_{k=1}^{n_{l+1}} \gamma_{l+1,k} q_{l+1,k}.
 		&=: -n_l b,
	\end{align*}
	and
	\[
		c - n_{l-1}c_{l-1} - n_{l+1}c_{l+1}
		= 1 - n_l + b.
	\]
	Moreover, the matrix $(\partial F_{l,k} / \partial q_{l,j})_{1 \le j,k\le
	n_l}$ is nonsingular as long as $b$ is not a non-positive integer bigger than
	$n_l$.
\end{lemma}

Indeed, a root of $P(z)$ is simple if and only if it does not coincide with $0$
or any $q_{l\pm 1, k}$.  If the roots $(q_{l,k})$ of $P(z)$ are all simple,
then they solve the equations~\cite{marden1966}
\begin{multline*}
	\sum_{1 \le k \ne j \le n_l} \frac{2}{q_{l,k}-q_{l,j}}
	+\sum_{1 \le j \le n_{l+1}} \frac{-c_{l+1}}{q_{l,k} - q_{l+1,j}}\\
	+\sum_{1 \le j \le n_{l-1}} \frac{-c_{l-1}}{q_{l,k} - q_{l-1,j}}
	+\frac{c}{q_{l,k}}
	= \frac{F_{l,k}}{q_{l,k}} = 0,
\end{multline*}
which is exactly our balance condition; see Remark~\ref{rem:electrob}.
Moreover, an equation system generalizing~\eqref{eq:hypergeometriccoeff} has
been obtained in~\cite[\S 136]{heine1878}, from which we may conclude the
non-singularity of the Jacobian.  In fact, there are
\[
	\binom{n_{l-1} + n_l + n_{l+1} - 1}{n_{l-1} + n_{l+1} - 1}
\]
choices of $\gamma$ for which~\eqref{eq:lame} has a polynomial solution of
degree $n_l$~\cite[\S 135]{heine1878}.

This observation allows us to easily construct balanced and rigid
configurations of type $(1,n,1)$:  Up to reparametrizations and complex
scalings, we may assume that $c_2 = 1$ and $q_{1,1} = 1$.  Then $q_{3,1}$ must
be real, and $(q_{2,k})_{1 \le k \le n}$ are given by roots of a Heun
polynomial.  Such a configuration depends locally on four real parameters,
namely $q_{3,1}$, $c_1$, $c_3$ and $c$ (or $b$).  When these are given, we have
$n+1$ Heun polynomials, each of which gives balanced positions of $q_{2,k}$'s.
For each of the Heun polynomials $P$, we have
\begin{align*}
	\dot\theta_{2,0} - \dot\theta_{1,0} &= \frac{P'(1)}{P(1)} - c_1,\\
	\dot\theta_{3,0} - \dot\theta_{2,0} &= c - 1 = b + c_1 + c_3 - n,\\
	\dot\theta_{4,0} - \dot\theta_{3,0} &= \frac{P'(q_{3,1})}{P(q_{3,1})} - c_3.
\end{align*}
Together with the family parameter $\tau$, the surface depends locally on five
parameters, which is expected because there are eight ends.

\begin{example}[Symmetric examples]
	When $q_{3,1} = q_{1,1} = 1$, the Heun polynomial reduces to a hypergeometric
	polynomial ${_2F_1}(-n, b; c; z)$, where $c_1 + c_3 = n-1-b+c$.  Assume
	further that $b+c=1-n$, so $c_1+c_3=-2b$.  This imposes a symmetry in the
	configuration.  Because $(c_1+c_3)P'(1)/P(1) = -nb$, the embeddedness
	conditions simplify to  
% The embeddedness conditions \eqref{eq:embedcond} are
% \begin{align*}
% 	\dot\theta_{1,0} > \dot\theta_{2,0} & \Rightarrow c_1(c_1+c_3)>-nb, \\
% 	\dot\theta_{1,\infty} > \dot\theta_{2,\infty} & \Rightarrow c_1(c_1+c_3)>n(c_1+c_3+b), \\
% 	\dot\theta_{2,0} > \dot\theta_{3,0} & \Rightarrow c<1, \\
% 	\dot\theta_{2,\infty} > \dot\theta_{3,\infty} & \Rightarrow b>-n,\\
% 	\dot\theta_{3,0} > \dot\theta_{4,0} & \Rightarrow c_3(c_1+c_3)>-nb, \\
% 	\dot\theta_{3,\infty} > \dot\theta_{4,\infty} & \Rightarrow c_3(c_1+c_3)>n(c_1+c_3+b).
% \end{align*}
% \begin{align*}
% 	\frac{n}{2} & < c_1 < \frac{3n}{2}+2c-2, \\
% 	c & < 1, \\
% 	b & > -n.
% \end{align*}
	\[
		c_1 > n/2,\quad c_3 > n/2,\quad 1-n/2 < c < 1,\quad -n < b < -n/2.
	\]
	As explained in Example \ref{ex:2}, the hypergeometric polynomial has real
	roots if $b$ and $c$ lie in the blue regions of Figure~\ref{fig:bcplot}.
	More specifically:
	\begin{itemize}
		\item When $n \ge 2$ and $0 < c < 1$, ${_2F_1}(-n, b; c; z)$ has $n$
			simple negative roots.  See Figure \ref{fig:151} for an example with $n=5$.

		\item When $n\ge 4$ and $-1<c<0$, or $n = 3$ and $-1/2 < c < 0$,
			${_2F_1}(-n, b; c; z)$ has $n-2$ simple negative roots, one root $0 < z_0
			< 1$, and another root $1/z_0 > 1$. \qed
	\end{itemize} 

	\begin{figure}[htbp] %%%%%%%%%%
		\centering
 		\includegraphics[height=2in]{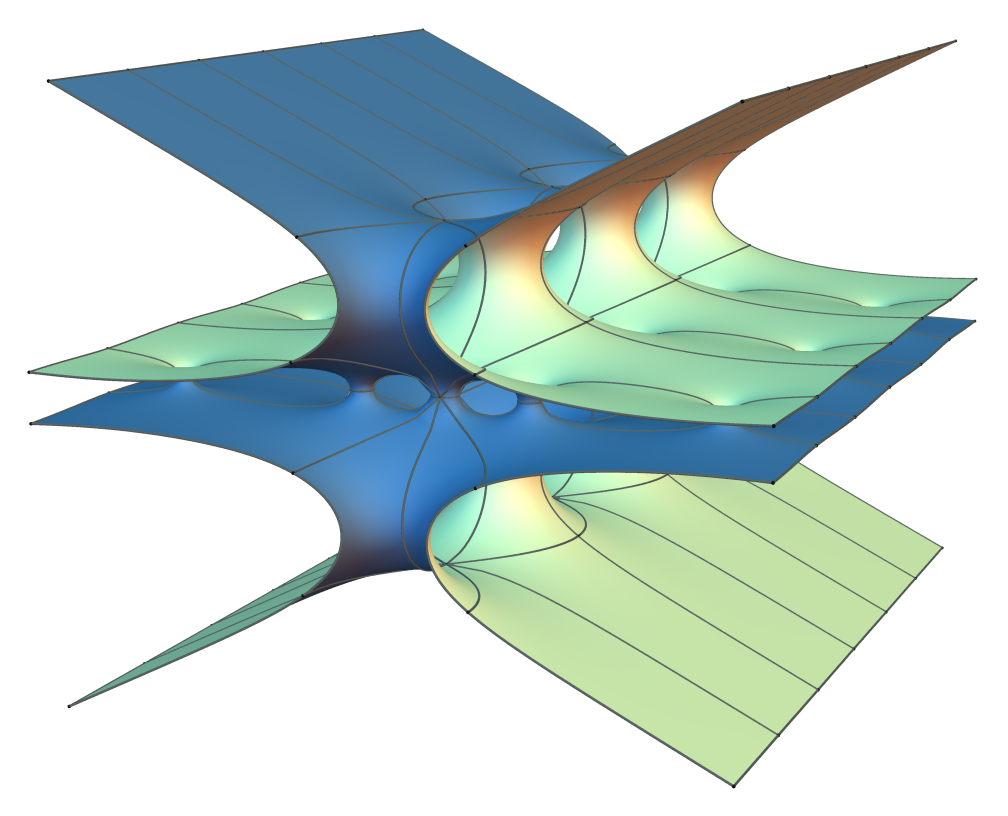}
		\caption{Genus 4 example with $n=5$ and $0<c<1$.}
		\label{fig:151}
	\end{figure} %%%%%%%%%%
\end{example}

\begin{example}[Offset handles]
	There are embedded examples in which the handles are not symmetrically
	placed.  For instance, one balanced configuration of type $(1,2,1)$ is given
	by
	\begin{align*}
		&q_{3,1} = 2/3,&& q_{2,1} = -1/3, && q_{2,2} = -23,\\
		&c_1 = 8/5, && c_3 = 4189/2890, && b = -9857/8670,
	\end{align*}
	so $c = 3956/4335$. \qed
	% $\dot\theta_{2,0}=\dot\theta_{1,0}-97/120$,
	% $\dot\theta_{3,0}=\dot\theta_{1,0}-2071/2312$, and
	% $\dot\theta_{4,0}=\dot\theta_{1,0}-4063763/2462280$.
\end{example}

\subsection{Concatenating surfaces of type (1,n,1)} \label{ssec:1n1n1}

We describe a family of examples in the same spirit
as~\cite[Proposition~2.3]{traizet2002b}.  Assume that we are in possession of
$R$ configurations of type $(1,n^{(r)},1)$, $n^{(r)} > 1$, $1 \le r \le R$.  In
the following, we use superscript $(r)$ to denote the parameters of the $r$-th
configuration.  Up to reparameterizations and complex scalings, we may assume
that $c_2^{(r)} = 1$ and $q_{1,1}^{(r)}=1$.  Then we may concatenate these
configurations into one of type
\[
	(1, n_2, 1, n_4, 1, \cdots, 1, n_{2R}, 1)
\]
such that $q_{1,1}=1$, $c_1 = 1$, and for $1 \le r \le R$, we have $n_{2r} =
n^{(r)}$,
\begin{align*}
	q_{2r,k} &= q_{2r-1,1} q_{2,k}^{(r)}, &
	q_{2r+1,1} &= q_{2r-1,1} q_{3,1}^{(r)},\\
	c_{2r} &= c_{2r-1}/c_1^{(r)}, &
	c_{2r+1} &= c_{2r-1} c_3^{(r)}/c_1^{(r)}, &
\end{align*}
and
\[
	\dot\theta_{2r+1, 0} - \dot\theta_{2r,0} = c_{2r} (c^{(r)} - 1).
\]
The balance of even layers then follows from the balance of each
sub-configuration.  The balance of odd layers leads to
\[ 
	\dot\theta_{2r, 0} - \dot\theta_{2r-1,0} =
	c_{2r} ( \dot\theta_{2,0}^{(r)} - \dot\theta_{1,0}^{(r)} + c_1^{(r)} )
	+ c_{2r-2} ( \dot\theta_{4,0}^{(r-1)} - \dot\theta_{3,0}^{(r-1)} + c_3^{(r-1)} )
	- c_{2r-1}
	% \frac{P'^{(r)}(1)}{P^{(r)}(1)}
	% + \frac{P'^{(r-1)}(q_{3,1}^{(r)})}{P^{(r-1)}(q_{3,1}^{(r)})} - c_{2r-1}
\]
for $1 \le r \le R+1$.  As expected, such a configuration depends locally on
$4R$ real parameters, namely $q_{3,1}^{(r)}$, $c_1^{(r)}$, $c_3^{(r)}$, and
$c^{(r)}$, $1 \le r \le R$.

We may impose symmetry by assuming that $q_{3,1}^{(r)} = 1$, so $q_{2r+1,1}=1$
for all $0 \le r \le R$, and that $b^{(r)}+c^{(r)} = 1 - n^{(r)}$, so
$c_1^{(r)} + c_3^{(r)} = n^{(r)} - 1 - b^{(r)} + c^{(r)} = -2b^{(r)}$.  Then
$q_{2r,k} = q_{2,k}^{(r)}$, $1 \le k \le n^{(r)}$, are given by the roots of
${_2F_1}(-n^{(r)}, b^{(r)}; c^{(r)}; z)$, $1 \le r \le R$.  Because
\[
 	\dot\theta_{2,0}^{(r)} - \dot\theta_{1,0}^{(r)} + c_1^{(r)} =
 	\dot\theta_{4,0}^{(r)} - \dot\theta_{3,0}^{(r)} + c_1^{(r)} = n^{(r)}/2,
\]
the embeddedness conditions simplify to
\[
	2 n_l c_l > n_{l-1} c_{l-1} + n_{l+1} c_{l+1}
\]
for $1 \le l \le L$ under the condition that $n_0 = n_{L+1} = 0$.  That is, the
sequence $(n_l c_l)_{1 \le l \le L}$ must be concave.  For even $l$, the
concavity implies that $b^{(r)} > -n$ hence $c^{(r)} < 1$ for all $1 \le r \le
R$.  We may choose, for instance, $n_l c_l = \ln(1+l)$ or $n_l c_l =
(\exp l - 1)/\exp(l-1)$ to obtain embedded minimal surfaces.

% \[
% 	b = n^{(r)}\Big(\frac{1}{4(4^r-1)}-1\Big),
% 	\quad \text{and} \quad
% 	c_1^{(r)} = n^{(r)} \frac{4^r-2}{4^r-1},
% \]
% so
% \[
% 	c = 1 - \frac{n^{(r)}}{4}\frac{1}{4^r-1}
% 	\quad \text{and} \quad
% 	c_3^{(r)} = n^{(r)} \frac{4^r-1/2}{4^r-1}.
% \]

\begin{remark}
	We can also append a configuration of type $(1,n^{(r)})$ to the sequence of
	$(1, n^{(r)}, 1)$-configurations to obtain a configuration of type
	\[
		(1,n_2,1,n_4,1,\ldots,1,n_{2R-2},1,n_{2R}),
	\]
	where the $q_{l,k}, c_l,\dot\theta_{l,0}$ terms are defined as above.
	Therefore, an embedded example of any genus with any even number ($>2$) of ends
	can be constructed.
\end{remark}

\subsection{Numerical examples}

The balance equations can be combined into one differential equation that is
much easier to solve.  A solution of this differential equation corresponds to
a lot of balance configurations that are equivalent by permuting the locations
of the nodes.  

\begin{lemma}
	Let $L$ be a positive integer, $n_1,n_2,\ldots,n_L\in\N$, and suppose $\{q_{l,k}\}$ is a configuration such that the $q_{l,k}$ are distinct.  Let
	\[
		P_l(z)=\prod_{k=1}^{n_l}(z-q_{l,k}),\quad P(z)=\prod_{l=1}^LP_l(z), \quad P_0(z)=P_{L+1}(z)=1,\\
	\]
	and
	\[
		\begin{split}
			\mathcal{F} P(z)=\sum_{l=1}^L&\left(\frac{c_l^2zP_l''(z)P(z)}{P_l(z)}-\frac{c_lc_{l+1}zP_l'(z)P_{l+1}'(z)P(z)}{P_l(z)P_{l+1}(z)}\right.\\
			&\hspace{.2in}\left.+\left(c_l^2+c_l(\dot\theta_{l+1,0}-\dot\theta_{l,0})\right)\frac{P_l'(z)P(z)}{P_l(z)}\right).
		\end{split}
	\]
	Then the configuration $\{q_{l,k}\}$ is balanced if and only if $\mathcal{F}P(z)\equiv 0$.
\end{lemma}

\begin{proof}
	We have seen that,
	\[
		\frac{P_l''(q_{l,k})}{P'_l(q_{l,k})}=\sum_{1\leq k\neq j\leq n_l}\frac{2}{q_{l,k}-q_{l,j}},
		\qquad
		\frac{P_{l\pm1}'(q_{l,k})}{P_{l\pm1}(q_{l,k})}=\sum_{j=1}^{n_{l\pm1}}\frac{1}{q_{l,k}-q_{l\pm1,j}}.
	\]
	Define
	\[
		F_l(z) = \frac{c_l^2 z P_l''(z)}{P_l'(z)}-\frac{c_lc_{l+1} z P_{l+1}'(z)}{P_{l+1}(z)}-\frac{c_lc_{l-1} z P_{l-1}'(z)}{P_{l-1}(z)}+c_l^2+c_l(\dot\theta_{l+1,0}-\dot\theta_{l,0}).
	\]
	Then $F_{l,k} = F_l(q_{l,k})$.
	% So we can express the force equation $F_{l,k}$ as
	% \[
	% 	\frac{c_l^2q_{l,k}P_l''(q_{l,k})}{P_l'(q_{l,k})}-\frac{c_lc_{l+1}q_{l,k}P_{l+1}'(q_{l,k})}{P_{l+1}(q_{l,k})}-\frac{c_lc_{l-1}q_{l,k}P_{l-1}'(q_{l,k})}{P_{l-1}(q_{l,k})}+c_l^2+c_l(\dot\theta_{l+1,0}-\dot\theta_{l,o}).
	% \]
	Set
	\[
		\begin{split}
			Q_l(z)&=\frac{P_l'(z)P(z)}{P_l(z)}F_l(z)\\
			&=\frac{c_l^2zP_l''(z)P(z)}{P_l(z)}-\frac{c_lc_{l+1}zP_l'(z)P_{l+1}'(z)P(z)}{P_l(z)P_{l+1}(z)}-\frac{c_lc_{l-1}zP_{l-1}'(z)P_l'(z)P(z)}{P_{l-1}(z)P_l(z)}\\
			&\hspace{.2in}+\left(c_l^2+c_l(\dot\theta_{l+1,0}-\dot\theta_{l,0})\right)\frac{P_l'(z)P(z)}{P_l(z)}.
		\end{split}
	\]
	Then $F_{l,k}=0$ if and only if $Q_l(q_{l,k})=0$.  

	Now observe that $Q_l(z)$ and $Q(z) = \mathcal{F} P(z)$ are polynomials with
	degree strictly less than
	\[
		\deg P = N = \sum_{l=1}^Ln_l,
	\]
	and $Q(q_{l,k})=Q_l(q_{l,k})$ for $1 \le k \le n_l$ and $1 \le l \le L$.  If
	$Q\equiv 0$ then $Q_l(q_{l,k})=0$ and so $\{q_{l,k}\}$ is a balanced
	configuration.  If $\{q_{l,k}\}$ is a balanced configuration then
	$Q(q_{l,k})=Q_l(q_{l,k})=F_{l,k}=0$.  Hence, $Q$ has at least $N$ distinct
	roots.  Since the degree of $Q$ is strictly less than $N$, we must have
	$Q\equiv 0$.
\end{proof}

It is relatively easy to numerically solve $\mathcal{F} P(z)\equiv 0$ as long
as we don't have too many levels and necks.  So we use this lemma to find
balanced configurations.  Since all previous examples admit a horizontal
reflection symmetry, we are most interested examples without this symmetry, or
with no non-trivial symmetry at all.

\begin{figure}[h]
  \centering
  \begin{subfigure}[b]{0.6\textwidth}
    \centering
    \includegraphics[width=\textwidth]{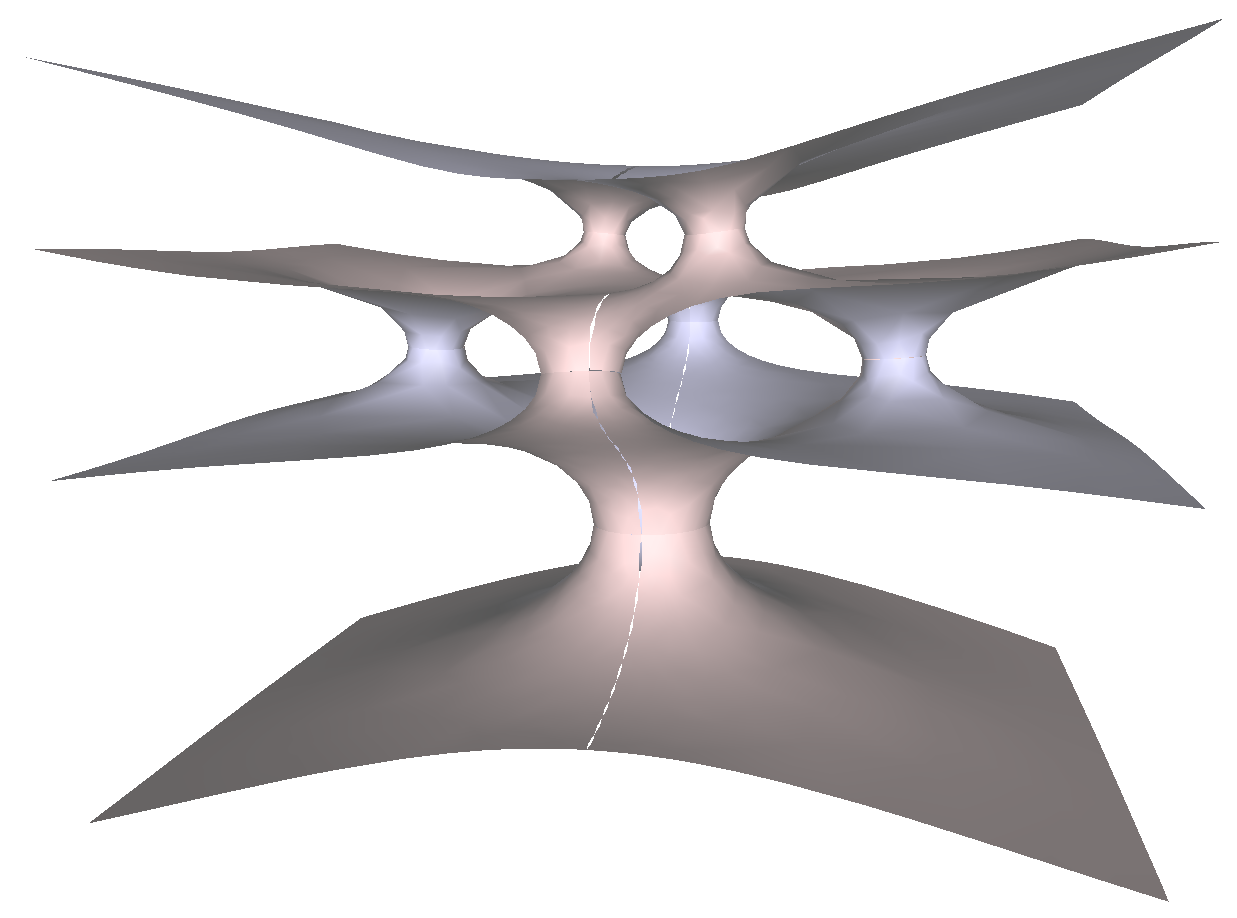}
  \end{subfigure}
  ~ %add desired spacing between images, e. g. ~, \quad, \qquad, \hfill etc. 
  \begin{subfigure}[b]{0.4\textwidth}
    \centering
    \includegraphics[width=\textwidth]{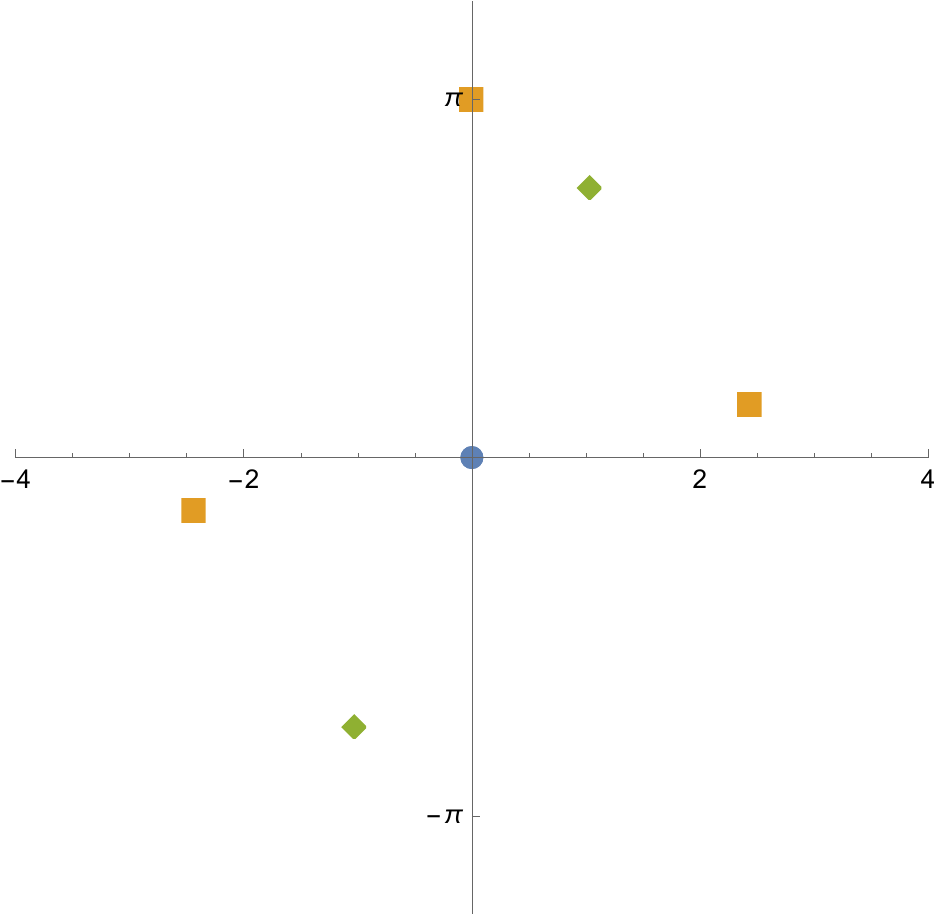}
  \end{subfigure}
  \caption{
  	A $(1,3,2)$ balanced configuration with no horizontal reflectional
  	symmetry.  The circles, squares, and diamonds represent the necks at levels
  	one, two, and three respectively.
	}\label{fig:132}
\end{figure}

Figure \ref{fig:132} shows an example with $L=3$,
\begin{gather*}
	n_1=1, n_2=3, n_3=2\\
	c_1=2, c_2=1, c_3=13/16\\
	\theta_{1,0}=0, \theta_{2,0}=-1/2, \theta_{3,0}=-27/16, \theta_{4,0}=-29/16.
\end{gather*}
This configuration corresponds to an embedded minimal surface with eight ends and
genus three in the quotient. It has no horizontal reflectional symmetry, but
does have a rotational symmetry.

\begin{figure}[h]
  \centering
  \begin{subfigure}[b]{0.5\textwidth}
    \centering
    \includegraphics[width=\textwidth]{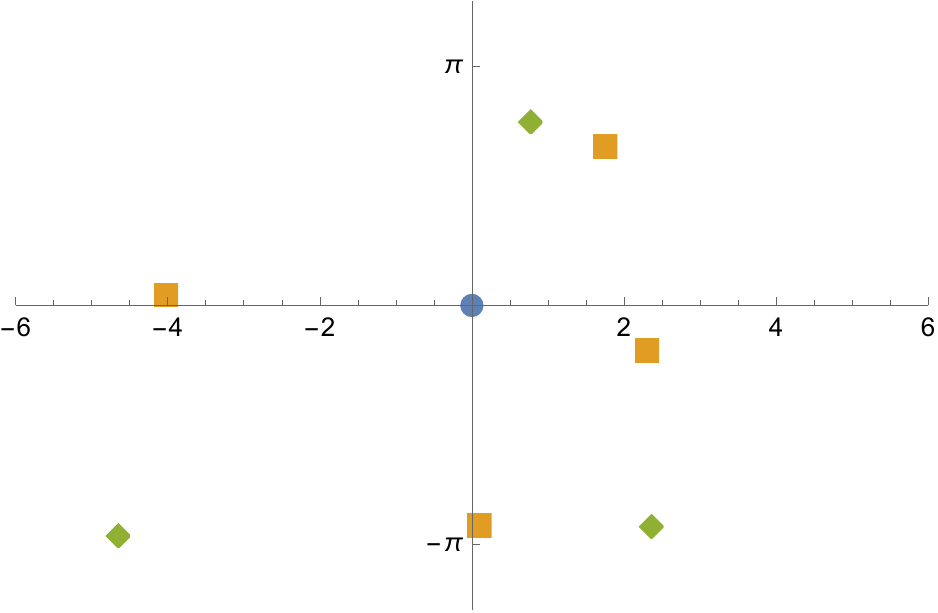}
  \end{subfigure}
  ~ %add desired spacing between images, e. g. ~, \quad, \qquad, \hfill etc. 
  \begin{subfigure}[b]{0.5\textwidth}
    \centering
    \includegraphics[width=\textwidth]{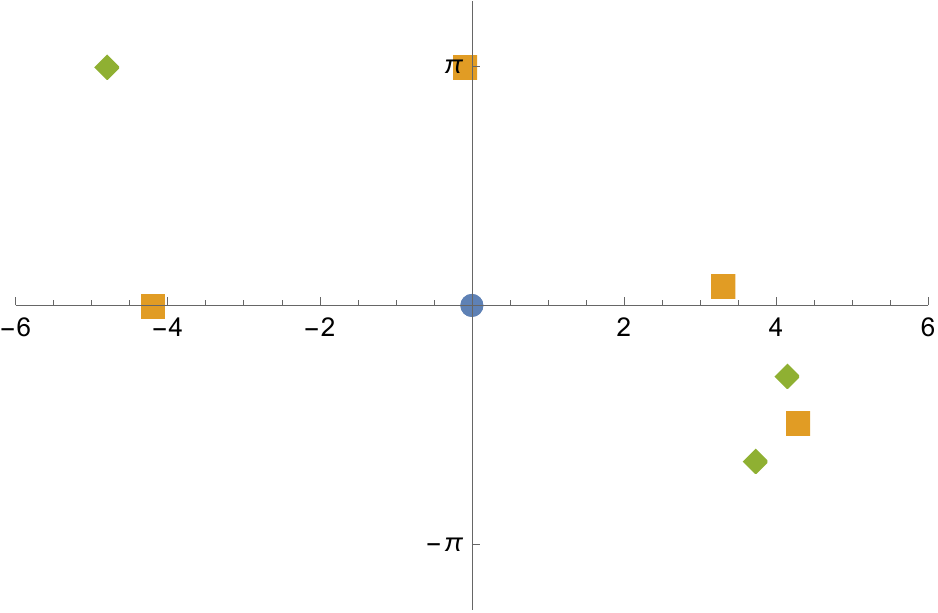}
  \end{subfigure}
  ~ %add desired spacing between images, e. g. ~, \quad, \qquad, \hfill etc. 
  \caption{
  	$(1,4,3)$ balanced configurations with no symmetries.  The circles, squares,
  	and diamonds represent the necks at levels one, two, and three respectively.
	}\label{fig:143}
\end{figure}

Figure \ref{fig:143} shows two examples with $L=3$,
\begin{gather*}
	n_1=1, n_2=4, n_3=3\\
	c_1=7/2, c_2=1, c_3=3/4\\
	\theta_{1,0}=0, \theta_{2,0}=-2, \theta_{3,0}=-13/5, \theta_{4,0}=-541/180.
\end{gather*}
These configurations correspond to embedded minimal surfaces with eight ends
and genus five in the quotient, with no nontrivial symmetry.

\begin{figure}[h]
  \centering
  \begin{subfigure}[b]{0.3\textwidth}
    \centering
    \includegraphics[width=\textwidth]{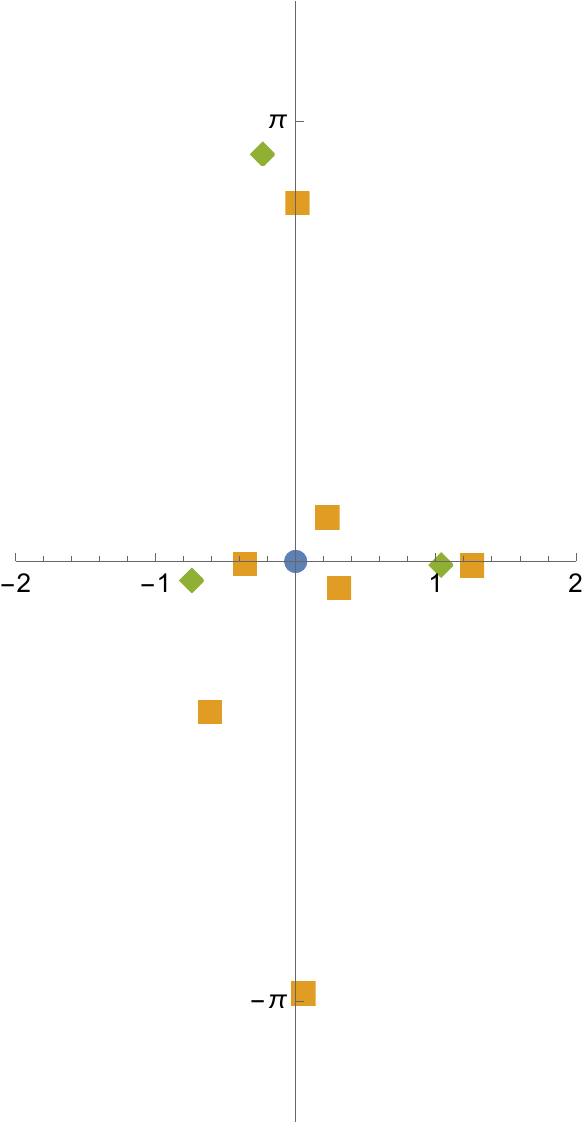}
  \end{subfigure}
  ~ %add desired spacing between images, e. g. ~, \quad, \qquad, \hfill etc. 
  \begin{subfigure}[b]{0.3\textwidth}
    \centering
    \includegraphics[width=\textwidth]{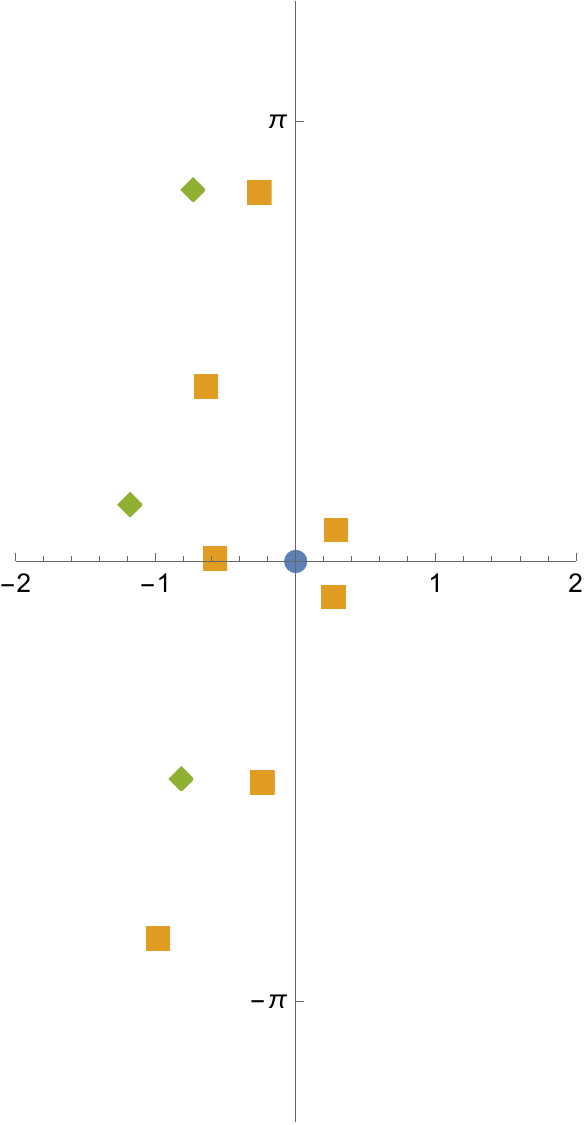}
  \end{subfigure}
  ~ %add desired spacing between images, e. g. ~, \quad, \qquad, \hfill etc. 
  \caption{
  	$(1,7,3)$ balanced configurations with no symmetries.  The circles, squares,
  	and diamonds represent the necks at levels one, two, and three respectively.
	}\label{fig:173}
\end{figure}

Figure \ref{fig:173} shows two examples with $L=3$,
\begin{gather*}
	n_1=1, n_2=7, n_3=3,\\
	c_1=17/7, c_2=1, c_3=3/2,\\
	\theta_{1,0}=0, \theta_{2,0}=-1/2, \theta_{3,0}=-3/2, \theta_{4,0}=-2468/441.
\end{gather*}
These configurations correspond to embedded minimal surfaces with eight ends
and genus eight in the quotient, with no nontrivial symmetry.

\section{Construction}

\subsection{Opening nodes}

To each vertical plane is associated a punctured complex plane $\C^\times_l
\simeq \C \setminus \{0\}$, $1 \le l \le L+1$.  They can be seen as Riemann
spheres $\hat\C_l \simeq \C \cup \{\infty\}$ with two fixed punctures at
$p_{l,0} = 0$ and $p_{l,\infty} = \infty$, corresponding to the two ends.

To each neck is associated a puncture $\cent p_{l,k} \in \C^\times_l$ and a
puncture $\cent{p'}_{l,k} \in \C^\times_{l+1}$.  Our initial surface at
$\tau=0$ is the noded Riemann surface $\Sigma_0$ obtained by identifying
$\cent p_{l,k}$ and $\cent{p'}_{l,k}$ for $1 \le l \le L$ and $1 \le k \le
n_l$.

As $\tau$ increases, we open the nodes into necks as follows:
Fix local coordinates $w_{l,0} = z$ in the neighborhood of $0 \in \hat\C_l$
and $w_{l,\infty} = 1/z$ in the neighborhood of $\infty \in \hat\C_l$.
For each neck, we consider parameters $(p_{l,k}, p'_{l,k})$ in the
neighborhoods of $(\cent p_{l,k}, \cent{p'}_{l,k})$ and local coordinates
\[
	% w_{l,k} = 2\frac{z-p_{l,k}}{z+p_{l,k}}
	w_{l,k} = \ln(z/p_{l,k})
	\quad\text{and}\quad
	% w'_{l,k} = 2\frac{z-p'_{l,k}}{z+p'_{l,k}}
	w'_{l,k} = \ln(z/p'_{l,k})
\]
in a neighborhood of $p_{l,k}$ and $p'_{l,k}$, respectively.  In this paper,
the branch cut of $\ln(z)$ is along the negative real axis, and we use the
principal value of $\ln(z)$ with imaginary part in the interval $(-\pi, \pi]$.

As we only open finitely many necks, we may choose $\delta > 0$ independent of
$k$ and $l$ such that the disks
\begin{gather*}
	|w_h| < 2\delta,\qquad h \in \gH \quad(=[1,L+1] \times \{0, \infty\}),\\
	|w_{l,k}| < 2\delta\quad\text{and}\quad|w'_{l,k}| < 2\delta,\qquad 1 \le l \le L, 1 \le k \le n_k
\end{gather*}
are all disjoint.  For parameters $t = (t_{l,k})_{1 \le l \le L, 1 \le k \le
n_l}$ in a neighborhood of $0$ with $|t_{l,k}| < \delta^2$, we remove the disks
\[
	|w_{l,k}|<|t_{l,k}|/\delta
	\quad\text{and}\quad
	|w'_{l,k}| < |t_{l,k}|/\delta
\]
and identify the annuli
\[
	|t_{l,k}|/\delta \le |w_{l,k}| \le \delta
	\quad\text{and}\quad
	|t_{l,k}|/\delta \le |w'_{l,k}| \le \delta
\]
by
\[
	w_{l,k} w'_{l,k} = t_{l,k}.
\]
If $t_{l,k} \ne 0$ for all $1 \le l \le L$ and $1 \le k \le n_l$, we obtain a
Riemann surface denoted by $\Sigma_t$.

\subsection{Weierstrass data}

We construct a conformal minimal immersion using the Weierstrass
parameterization in the form
\[
	z\mapsto	\re \int^z (\Phi_1, \Phi_2, \Phi_3) ,
\]
where $\Phi_i$ are meromorphic 1-forms on $\Sigma_t$ satisfying the
conformality equation
\begin{equation}\label{eq:conformal}
	Q := \Phi_1^2 + \Phi_2^2 + \Phi_3^2 = 0.
\end{equation}

\subsubsection{A-periods}

We consider the following fixed domains in all $\Sigma_t$:
\begin{align*}
	U_{l,\delta}=\{ z \in \hat\C_v \colon
		&|\cent w_{l,k}(z)| > \delta/2 \quad \forall 1 \le k \le n_l \quad \text{if}\; 1 \le l \le L \\
	\text{and}\quad &|\cent{w'}_{l,k}(z)| > \delta/2 \quad \forall 1 \le k \le n_{l-1} \quad \text{if}\; 2 \le l \le L+1\}
\end{align*}
and $U_{\delta}=\bigsqcup_{1 \le l \le L} U_{l,\delta}$.

Let $A_{l,k}$ denote a small counterclockwise circle in $U_{l,\delta}$ around
$p_{l,k}$; it is then homologous in $\Sigma_t$ to a clockwise circle in
$U_{l+1,\delta}$ around $p'_{l,k}$.  Moreover, let $A_{l,0}$ (resp.\
$A_{l,\infty}$) denote a small counterclockwise circle in $U_{l,\delta}$ around
$0$ (resp.\ $\infty$).

Recall that the vertical period vector is assumed to be $(0,0,2\pi)$, so we
need to solve the A-period problems
\[
	\re \int_{A_h} (\Phi_1, \Phi_2, \Phi_3) = (0,0, 2\pi\sigma_h)
	\quad\text{and}\quad
	\re \int_{A_{l,k}} (\Phi_1, \Phi_2, \Phi_3) = (0,0,0)
\]
for $h \in \gH$, $1 \le l \le L$, and $1 \le k \le n_l$.  Here, the orientation
$\sigma_h = \pm 1$ satisfies
\[
	\sigma_h = - \sigma_{\rotate(h)},
\]
where the ``counterclockwise rotation'' $\rotate$ on $\gH$ is defined by
\begin{equation}\label{eq:degrotate}
	\begin{cases}
		\rotate(0_l) = 0_{l-1}, & 2 \le l \le L+1,\\
		\rotate(0_1)= \infty_1,\\
		\rotate(\infty_l) = \infty_{l+1}, & 1 \le l \le L,\\
		\rotate(\infty_{L+1}) = 0_{L+1}.
	\end{cases}
\end{equation}
In particular, we have $\sigma_{l,0} = - \sigma_{l,\infty}$ for all $1 \le l
\le L+1$.

Recall that the surface tends to an $(L+1)$-sheeted $xz$-plane in the limit
$\tau \to 0$.  So we define the meromorphic functions $\Phi_1$, $\Phi_2$ and
$\Phi_3$ as the unique regular 1-forms on $\Sigma_t$ (see~\cite[\S
8]{traizet2013}) with simple poles at $p_h$, $h \in \gH$, and the A-periods
\begin{align*}
	\int_{A_h}{(\Phi_1, \widetilde\Phi_2, \Phi_3)} &=
	2 \pi \ii ( \alpha_h, \beta_h, \gamma_h - \ii \sigma_h), && h \in \gH, \\
	\int_{A_{l,k}} (\Phi_1, \widetilde\Phi_2, \Phi_3) &=
	2 \pi \ii (\alpha_{l,k}, \beta_{l,k}, \gamma_{l,k}), && 1 \le l \le L, 1 \le k \le n_l,
\end{align*}
where $\Phi_2 = \tau\widetilde\Phi_2$ and, by Residue Theorem, it is necessary
that
\begin{align}
	\alpha_{l,0} + \alpha_{l,\infty} + \sum_{1 \le k \le n_l} \alpha_{l,k} -
	\sum_{1 \le k \le n_{l-1}} \alpha_{l-1,k}&= 0\label{eq:resalpha},\\
	\beta_{l,0} + \beta_{l,\infty} + \sum_{1 \le k \le n_l} \beta_{l,k} -
	\sum_{1 \le k \le n_{l-1}} \beta_{l-1,k}&= 0\label{eq:resbeta},\\
	\gamma_{l,0} + \gamma_{l,\infty} + \sum_{1 \le k \le n_l} \gamma_{l,k} -
	\sum_{1 \le k \le n_{l-1}} \gamma_{l-1,k}&= 0\label{eq:resgamma},
\end{align}
for $1 \le l \le L+1$.  Then the A-period problems are solved by definition.

\subsubsection{Balance of ends}

Summing up~\eqref{eq:resbeta} over $l$ gives
\begin{equation}\label{eq:endbalanceb}
	\sum_{h \in \gH} \beta_h = 0,
\end{equation}
which we use to replace Equation \eqref{eq:resbeta} with $l = L+1$.

In this paper, the punctures $p_{l,0}$ and $p_{l,\infty}$ correspond to Scherk-type
ends.  Hence we fix
\begin{equation}\label{eq:fixalphagamma}
	\alpha^2_h + \tau^2 \beta_h^2 \equiv 1
	\quad\text{and}\quad
	\gamma_h \equiv 0
\end{equation}
for all $h \in \gH$, so that (the stereographic projection of) the Gauss map $G
= -(\Phi_1 + \ii \Phi_2)/\Phi_3$ extends holomorphically to the punctures $p_h$
with unitary values.  Then Equations \eqref{eq:resgamma} are not independent:
if it is solved for $1 \le l \le L$, it is automatically solved for $l = L+1$.  

In particular, at $\tau=0$, we have $\alpha_h^2=1$.  In view of the orientation
of the ends, we choose $\alpha_{l,0} = 1$ and $\alpha_{l,\infty} = -1$
so that $G(p_{l,\infty}) = G(p_{l,0}) = \ii\sigma_{l,0}$.

Summing up~\eqref{eq:resalpha} over $l$ gives
\begin{equation}\label{eq:endbalancea}
	\sum_{1 \le l \le L+1} \Big( \sqrt{1-\tau^2\beta_{l,\infty}^2} - \sqrt{1-\tau^2\beta_{l,0}^2} \Big) = 0,
\end{equation}
which we use to replace Equation~\eqref{eq:resalpha} with $l = L+1$.
\begin{remark}
	The conditions~\eqref{eq:endbalanceb} and~\eqref{eq:endbalancea} are
	disguises of the balance condition of Scherk ends, namely that the unit
	vectors in their directions should sum up to $0$.
\end{remark}

\subsubsection{B-periods}

For $1 \le l \le L+1$, we fix a point $O_l \in U_{l,\delta}$.  For every $1 \le
l \le L$ and $1 \le k \le n_l$ and $t_{l,k}\neq 0$, let $B_{l,k}$ be the
concatenation of
\begin{enumerate}
	\item a path in $U_{l,\delta}$ from $O_l$ to $w_{l,k} = \delta$,
	\item the path parameterized by $w_{l,k} = \delta^{1-2s}\,t_{l,k}^s$ for $s \in
		[0,1]$, from $w_{l,k} = \delta$ to $w_{l,k} = t_h/\delta$, which is identified with
		$w'_{l,k}=\delta$, and
	\item a path in $U_{l+1,\delta}$ from $w'_{l,k}=\delta$ to $O_{l+1}$.
\end{enumerate}

We need to solve the B-period problem, namely that
\begin{equation}\label{eq:Bperiod}
	\re \int_{B_{l,k}} (\Phi_1, \Phi_2, \Phi_3) =
	\re \int_{B_{l,1}} (\Phi_1, \Phi_2, \Phi_3).
\end{equation}

\subsubsection{Conformality}

\begin{lemma}
	For $t$ sufficiently close to $0$, the conformality condition
	\eqref{eq:conformal} is equivalent to
	\begin{align}
		\fG_{l,k} := \int_{A_{l,k}} \frac{w_{l,k} Q}{dw_{l,k}} &= 0, \quad 1 \le l \le L, \quad 1 \le k \le n_l, \label{eq:conform1}\\
		\fF_{l,k} := \int_{A_{l,k}} \frac{Q}{dw_{l,k}} &= 0, \quad 1 \le l \le L,\quad 2 \le k \le n_l, \label{eq:conform2}\\
		\fF'_{l,k} := \int_{A'_{l,k}} \frac{Q}{dw'_{l,k}} &= 0, \quad 1 \le l \le L,\quad 1 + \delta_{l,L} \le k \le n_l, \label{eq:conform3}
	\end{align}
	where $A'_{l,k}$ in~\eqref{eq:conform3} denotes a small counterclockwise
	circle in $U_{l+1,\delta}$ around $p'_{l,k}$ (hence homologous to
	$-A_{l,k}$), and $\delta_{l,L}=1$ if $l = L$ and $0$ otherwise.
\end{lemma}

\begin{proof}
	By our choice of $\alpha_h$ and $\gamma_h$, the quadratic differential $Q$
	has at most simple poles at the $2L+2$ punctures $p_h$, $h \in \gH$.  The
	space of such quadratic differentials is of complex dimension
	$3(N-L)-3+(2L+2) = 3N-L-1$.  We will prove that
	\[
		Q \mapsto (\fG, \fF, \fF')
	\]
	is an isomorphism.  We prove the claim at $t=0$; then the claim follows by
	continuity.

	Consider $Q$ in the kernel.  Recall from~\cite{traizet2008} that a regular
	quadratic differential on $\Sigma_0$ has at most double poles at the nodes
	$p_{l,k}$ and $p'_{l,k}$.  Then~\eqref{eq:conform1} guarantees that $Q$ has
	at most simple poles at the nodes.  By~\eqref{eq:conform2}
	and~\eqref{eq:conform3}, $Q$ may only have simple poles at $p_{l,1} \in
	\C^\times_l$, $1 \le l \le L$, and $p'_{L,1} \in \C^\times_{L+1}$.  So, on
	each Riemann sphere $\hat\C_l$, $Q$ is a quadratic differential with at most
	simple poles at three punctures; the other two being $0$, $\infty$.  But such
	a quadratic differential must be $0$.
\end{proof}

% \subsubsection{Parameter count}
%
% {\color{blue}
% 	First count the parameters
%
% 	\begin{itemize}
% 		\item We have a parameter $\tau$.
%
% 		\item We have $2N$ punctures $p$ and $p'$.  Up to complex scaling, these are
% 			$2N-L-1$ complex numbers.
%
% 		\item We have $N$ complex numbers $t$.
%
% 		\item We have the residues $\alpha$, $\beta$, and $\gamma$.  The $\alpha_h$ and
% 			$\gamma_h$ has been fixed, only $\beta_h$ are free.  So these are
% 			$3N+2(L+1) = 3N+2L+2$ real parameters.
% 	\end{itemize}
% 	So there are $2(3N-L-1)+3N+2L+2+1 = 9N+1$ real parameters.
%
% 	Then count the equations
%
% 	\begin{itemize}
% 		\item The conformality equations $3N-L-1$ complex equations.
%
% 		\item The B-periods give $3N-3L$ real equations.
%
% 		\item $3L+2$ Residue Theorems.
% 	\end{itemize}
%
% 	So there are $2(3N-L-1)+3N+2 = 9N-2L$ real equations.
%
% 	There are $2L+1$ more parameters than equations. It should however be
% 	$2(L+1)-3 = 2L-1$, so we have two more.  One comes from the rotation.
% 	Another is from the fact that $(\beta,\tau) \to
% 	(\lambda\beta,\tau/\lambda)$ gives the same Weierstrass data.
% }

\subsection{Using the Implicit Function Theorem}

All parameters varies in a neighborhood of their central values, denoted by a
superscript $\circ$.  We will see that
\[
	\cent\beta_h = \dot\theta_h,\;
	\cent\alpha_{l,k}=\cent\gamma_{l,k} = 0,\;
	\cent\beta_{l,k} = -c_l,\;
	\cent{p'}_{l,k} = \overline{p_{l,k}}.
\]

Let us first solve Equations~\eqref{eq:endbalanceb} and~\eqref{eq:endbalancea}.
\begin{proposition}\label{prop:theta}
	Given a configuration $(q, \dot\theta)$ such that $\Theta_1=\Theta_2=0$.  For
	$\tau$ sufficiently small and $\beta_h$ close to $\cent\beta_h =
	\dot\theta_h$, the solutions $(\tau,\beta)$ to
	Equations~\eqref{eq:endbalanceb} and~\eqref{eq:endbalancea} form a smooth
	manifold of dimension $2L+1$.
\end{proposition}

\begin{proof}
	At $\tau=0$, Equation~\eqref{eq:endbalanceb} is solved by $\cent\beta_h =
	\dot\theta_h$ if $\Theta_1=0$.  Taking the derivative
	of~\eqref{eq:endbalancea} with respect to $\tau^2$ gives
	\begin{equation}\label{eq:endbalancea0}
		\sum_{1 \le l \le L+1} \frac{\beta_{l,\infty}^2 - \beta_{l,0}^2}{2}  = 0,
	\end{equation}
	which is solved by $\cent\beta_h = \dot\theta_h$ if $\Theta_2=0$.  The
	proposition then follows from Lemma~\ref{lem:rigid} and the Implicit Function
	Theorem.
\end{proof}

From now on, we assume that the parameters $(\tau,(\beta_h)_{h\in\gH})$ are
solutions to Equations~\eqref{eq:endbalanceb} and~\eqref{eq:endbalancea} in a
neighborhood of $(0, \dot\theta)$.

\subsubsection{Solving conformality problems}

\begin{proposition}\label{prop:talphagamma}
	For $\tau$ sufficiently small and $\beta_{l,k}$, $p_{l,k}$, and $p'_{l,k}$ in
	a neighborhood of their central values, there exist unique values of
	$t_{l,k}$, $\alpha_{l,k}$, and $\gamma_{l,k}$, depending real-analytically on
	$(\tau^2, \beta, p, p')$, such that the balance equations~\eqref{eq:resalpha}
	and~\eqref{eq:resgamma} with $1 \le l \le L$ and the conformality Equations
	\eqref{eq:conform1} and~\eqref{eq:conform2} are solved.  Moreover, at
	$\tau=0$, we have $t_{l,k}=0$, $\alpha_{l,k} = \gamma_{l,k} = 0$,
	\[
		\frac{\partial t_{l,k}}{\partial(\tau^2)} = \frac{1}{4}\beta_{l,k}^2,
	\]
	and, for $2 \le k \le n_l$,
	\begin{equation}\label{eq:dalphagamma}
		\frac{\partial}{\partial(\tau^2)} (\alpha_{l,k} - \ii\sigma_{l,0} \gamma_{l,k})
		= -\frac{1}{2} \Res{\frac{\widetilde\Phi_2^2}{dw_{l,k}}}{p_{l,k}}
		= -\frac{1}{2} \Res{\frac{z\widetilde\Phi_2^2}{dz}}{p_{l,k}}.
	\end{equation}
\end{proposition}

Note that, according to this proposition, if $\cent\beta_{l,k} \ne 0$, then
$t_{l,k} > 0$ for sufficiently small $\tau$.

\begin{proof}
	At $\tau=0$, for $2 \le k \le n_l$ we have
	\[
		\fG_{l,k} = \int_{A_{l,k}} \frac{w_{l,k}Q}{dw_{l,k}} = 
		2\pi\ii(\alpha_{l,k}^2 + \gamma_{l,k}^2) = 0
	\]
	which vanishes when
	\[
		\alpha_{l,k} = \gamma_{l,k} = 0.
	\]
	Recall that $\alpha_h = \pm 1$ at $\tau=0$ and that $\gamma_h \equiv
	0$.  Then by Residue Theorem, we have
	\[
		\alpha_{l,1} = \gamma_{l,1} = 0.
	\]
	As a consequence, we have at $\tau=0$
	\[
		\cent\Phi_1 = dz/z,\quad
		\cent\Phi_2 = 0,\quad\text{and}\quad
		\cent\Phi_3 = - \ii \sigma_{l,0} dz/z,
	\]
	so $Q = 0$ as we expect.

	We then compute the partial derivatives at $\tau=0$
	\begin{align*}
		\frac{\partial}{\partial \alpha_{l,k}} \fF_{l,k}
		&=\int_{A_{l,k}} 2\frac{\cent\Phi_1}{dw_{l,k}}\frac{\partial\Phi_1}{\partial\alpha_{l,k}} \bigg|_{\tau=0}
		=\int_{A_{l,k}} 2\frac{dz/z}{dz/p}\frac{dz}{z-p_{l,k}}
		=4\pi\ii\,,\\
		\frac{\partial}{\partial \gamma_{l,k}} \fF_{l,k}
		&=\int_{A_{l,k}} 2\frac{\cent\Phi_3}{dw_{l,k}}\frac{\partial\Phi_3}{\partial\gamma_{l,k}} \bigg|_{\tau=0}
		=\int_{A_{l,k}} -2i\sigma_{l,0}\frac{dz/z}{dz/p}\frac{dz}{z-p_{l,k}}
		=4\pi\sigma_{l,0}\,,\\
		\frac{\partial}{\partial t_{l,k}} \fG_{l,k}
		&=\int_{A_{l,k}} \frac{2 w_{l,k}}{dw_{l,k}} \Big(
			\cent\Phi_1 \frac{\partial\cent\Phi_1}{\partial t_{l,k}}
			+ \cent\Phi_3 \frac{\partial\cent\Phi_3}{\partial t_{l,k}}
		\Big)\\
		&=\frac{-1}{\pi\ii} \Big(
			\int_{A_{l,k}} \frac{\cent\Phi_1}{w_{l,k}}
			\int_{A'_{l,k}} \frac{\cent\Phi_1}{w'_{l,k}} +
			\int_{A_{l,k}} \frac{\cent\Phi_3}{w_{l,k}}
			\int_{A'_{l,k}} \frac{\cent\Phi_3}{w'_{l,k}}
		\Big)\quad\text{\cite[Lemma~3]{traizet2008}}\\
		& = -8\pi\ii.
		% =\frac{-1}{\pi\ii}
		% [ (2\pi\ii)(2\pi\ii) + (2\pi\sigma_{l,0})(2\pi\sigma_{l+1,0}) ]
    %
		% \frac{w_{l,k}}{dw_{l,k}} \Big(
		% 	\Phi_1 \frac{\partial\Phi_1}{\partial t_{l,k}}
		% 	+ \Phi_3 \frac{\partial\Phi_3}{\partial t_{l,k}}
		% \Big)
    %
	\end{align*}
	All other partial derivatives vanish.  Therefore, by the Implicit Function
	Theorem, there exist unique values of $\alpha_{l,k}$, $\gamma_{l,k}$ (with $2
	\le k \le n_l$), and $t_{l,k}$ (with $1 \le k \le n_l$) that solve the
	conformality Equations \eqref{eq:conform1} and~\eqref{eq:conform2}.  Recall
	that $\alpha_h$ are determined by~\eqref{eq:fixalphagamma}. Then
	$\alpha_{l,1}$ and $\gamma_{l,1}$ are uniquely determined by the linear
	balance equations~\eqref{eq:resalpha} and~\eqref{eq:resgamma}.

	Moreover,
	\[
		\frac{\partial}{\partial(\tau^2)} \fF_{l,k} = \int_{A_{l,k}} \frac{\widetilde\Phi_2^2}{dw_{l,k}},\quad
		\frac{\partial}{\partial(\tau^2)} \fG_{l,k} = 2\pi\ii \beta_{l,k}^2.
	\]
	Hence the total derivatives
	\[
		\frac{d}{d(\tau^2)} \fF_{l,k} =
		4\pi\ii \frac{\partial\alpha_{l,k}}{\partial(\tau^2)} +
		4\pi\sigma_{l,0} \frac{\partial\gamma_{l,k}}{\partial(\tau^2)} +
		2\pi\ii \Res{\frac{\widetilde\Phi_2^2}{dw_{l,k}}}{p_{l,k}} = 0
	\]
	and
	\[
		\frac{d}{d(\tau^2)} \int_{A_{l,k}} \fG_{l,k} =
 		-8\pi\ii\frac{\partial t_{l,k}}{\partial(\tau^2)} + 2\pi\ii \beta_{l,k}^2 = 0.
 	\]
 	This proves the claimed partial derivatives with respect to $\tau^2$.
\end{proof}

\begin{remark}
	We see from the computations that our local coordinates $w$ and $w'$ are
	chosen for convenience.  Had we used other coordinates, the computations
	would be very difference, but $\partial(\alpha_{l,k} - \ii\sigma_{l,0}
	\gamma_{l,k})/\partial(\tau^2)$ would be invariant, and $\partial
	t_{l,k}/\partial(\tau^2)$ would be rescaled to keep the conformal type of
	$\Sigma_t$ (to the first order).  So the choice of local coordinates has no
	substantial impact on our construction.
\end{remark}

% From now on, we assume that the parameters $t_{l,k}$, $\alpha_{l,k}$, and
% $\gamma_{l,k}$ are given by Proposition~\ref{prop:talphagamma}.

\subsubsection{Solving B-period problems}

In the following, we make a change of variable $\tau = \exp(-1/\xi^2)$.

\begin{proposition}\label{prop:beta}
	Assume that the parameters $t_{l,k}$, $\alpha_{l,k}$, and $\gamma_{l,k}$ are
	given by Proposition~\ref{prop:talphagamma}.  For $\xi$ sufficiently small
	and $p_{l,k}$ and $p'_{l,k}$ in a neighborhood of their central values, there
	exist unique values of $\beta_{l,k}$, depending smoothly on $(\xi,p,p')$ and
	$(\beta_h)_{h \in \gH}$, such that the balance equation~\eqref{eq:resbeta}
	with $1 \le l \le L$ and the $y$-component of the B-period
	problem~\eqref{eq:Bperiod} are solved.  Moreover, at $\xi=0$ and $\beta_h =
	\cent\beta_h = \dot\theta_h$, we have $\beta_{l,k} = \beta_{l,1} = -c_l$
	where $c_l$ is given by~\eqref{eq:resc}.
\end{proposition}

\begin{proof}
	By Lemma 8.3 of~\cite{saddle1}, 
	\[
		\Big(\int_{B_{l,k}} \widetilde\Phi_2\Big)- \beta_{l,k} \ln t_{l,k}
	\]
	extends holomorphically to $t=0$ as bounded analytic functions of other
	parameters.  We have seen that $t_{l,k} \sim \tau^2\beta^2_{l,k}/4$.  So
	\[
		\mathfrak{H} := -\frac{\xi^2}{2} \re \Big( \int_{B_{l,k}}
		\widetilde\Phi_2 - \int_{B_{l,1}} \widetilde\Phi_2 \Big) = \beta_{l,k} - \beta_{l,1}
	\]
	at $\xi=0$.  Therefore, $\mathfrak{H}=0$ is solved at $\xi=0$ by $\beta_{l,k} =
	\beta_{l,1}$ for all $2 \le k \le n_l$, and $\beta_{l,1}= -c_l$ follows as
	\eqref{eq:resc} is just a reformulation of~\eqref{eq:resbeta}.  The
	proposition then follows by the Implicit Function Theorem.
\end{proof}

% From now on, we assume that the parameters $\beta_{l,k}$ are given by
% Proposition~\ref{prop:beta}.

\begin{proposition}\label{prop:pprime}
	Assume that the parameters $t_{l,k}$, $\alpha_{l,k}$, $\beta_{l,k}$ and
	$\gamma_{l,k}$ are given by Propositions~\ref{prop:talphagamma}
	and~\ref{prop:beta}.  For $\xi$ sufficiently small and $p_{l,k}$ in a
	neighborhood of their central values, there exist unique values of
	$p'_{l,k}$, depending smoothly on $\xi$, $p$, and $(\beta_h)_{h \in \gH}$,
	such that the $x$- and $z$-components of the B-period
	problem~\eqref{eq:Bperiod} are solved.  Moreover, up to complex scalings on
	$\C^\times_{l+1}$, $1 \le l \le L$, we have $p'_{l,k} = \overline{p_{l,k}}$
	at $\xi=0$ for any $1 < k \le n_l$.
\end{proposition}

\begin{proof}
	At $\xi=0$, recall that $\Phi_1 = dz/z$ and $\Phi_3 =
	-\ii\sigma_{l,0} dz/z$.  So
	\begin{align*}
		\re \int_{B_{l,k}} \Phi_1 - \re\int_{B_{l,1}} \Phi_1
		&= \re\ln(p_{l,k}/p_{l,1}) - \re\ln(p'_{l,k}/p'_{l,1}); \\
		\re \int_{B_{l,k}} \Phi_3 - \re \int_{B_{l,1}} \Phi_3
		&= \sigma_{l,0} (\im\ln(p_{l,k}/p_{l,1}) + \im\ln(p'_{l,k}/p'_{l,1})).
	\end{align*}
	They vanish if and only if $\ln(p_{l,k}/p_{l,1}) =
	\overline{\ln(p'_{l,k}/p'_{l,1})}$.  We normalize the complex scaling on
	$\C^\times_{l+1}$, $1 \le l \le L$, by fixing $p'_{l,1} =
	\overline{p_{l,1}}$. Then the B-period problem is solved at $\xi=0$ with
	$p'_{l,k} = \overline{p_{l,k}}$.  By the same argument as
	in~\cite{traizet2008}, the integrals are smooth functions of $\xi$ and other
	parameters, so the proposition follows by the Implicit Function Theorem.
\end{proof}

% From now on, we assume that the parameters $p'$ are given by
% Proposition~\ref{prop:beta}.

\subsubsection{Balancing conditions}

Define
\[
	\fR_{l,k} = \Res{\frac{z\widetilde\Phi_2^2}{dz}}{p_{l,k}}
	\quad\text{and}\quad
	\fR'_{l,k} = \Res{\frac{z\widetilde\Phi_2^2}{dz}}{p'_{l,k}}.
\]
Let the central values $\cent p_{l,k} = \conj^l q_{l,k}$, where $q$ is from a
balanced configuration.  So the central values $\cent{p'}_{l,k} = \conj^{l+1}
q_{l,k}$ and
\[
	\cent{\widetilde\Phi}_2 = \begin{cases}
		\overline{\conj^* \psi_l} & \text{on $\C^\times_l$ for $l$ odd,}\\
		\psi_l & \text{on $\C^\times_l$ for $l$ even.}
	\end{cases}
\]
Then we have
\[
	\overline{\fR_{l,k}} + \fR'_{l,k} = 2 \conj^{l+1} F_{l,k}
\]
at the central values, where $F_{l,k}$ is the force given by~\eqref{eq:force}.
Moreover, by Residue Theorem on $\C^\times_l$,
\begin{equation}\label{eq:residuephi2}
	\sum_{k=1}^{n_{l-1}}\fR'_{l-1,k} + \sum_{k=1}^{n_l}\fR_{l,k} +
	\beta_{l,0}^2 - \beta_{l,\infty}^2 = 0.
\end{equation}

\begin{proposition}\label{prop:extendF}
	Assume that the parameters $t_{l,k}$, $\alpha_{l,k}$, and $\gamma_{l,k}$ are
	given as analytic functions of $\tau^2$ by
	Proposition~\ref{prop:talphagamma}.  Then $\widetilde\fF'_{l,k} := \tau^{-2}
	\fF'_{l,k}$ extends analytically to $\tau=0$ with the value
	\[
		\begin{dcases}
 			4\pi\ii \conj^{l+1} F_{l,k}, & 2 \le k \le n_l,\\
			4\pi\ii \conj^{l+1} \Big(
				F_{l,1} + \sum_{j=1}^{l-1}\sum_{k=1}^{n_j}F_{j,k}
			\Big), & k = 1.
		\end{dcases}
	\]
\end{proposition}

\begin{proof}
	If $f(z)$ is an analytic function in $z$ and $f(0)=0$, then $f(z)/z$ extends
	analytically to $z=0$ with the value $df/dz\mid_{z=0}$.  We compute at $\tau
	= 0$ that
	\[
		\frac{\partial}{\partial \alpha} \fF'_{l,k} = -4\pi\ii
		\quad\text{and}\quad
		\frac{\partial}{\partial \gamma} \fF'_{l,k} = 4\pi\sigma_{l,0}.
	\]
	Then
	\[
		\frac{d}{d(\tau^2)} \fF'_{l,k}
		= -4\pi\ii \frac{\partial\alpha_{l,k}}{\partial(\tau^2)} +
		4\pi\sigma_{l,0} \frac{\partial\gamma_{l,k}}{\partial(\tau^2)} +
		2\pi\ii \fR'_{l,k}.
	\]
	For $2 \le k \le n_l$, by~\eqref{eq:dalphagamma}, $\widetilde\fF'_{l,k} :=
	\tau^{-2} \fF'_{l,k}$ extends to $\tau=0$ with the value
	\[
		\frac{d}{d(\tau^2)} \fF'_{l,k}
	 	= 2\pi\ii (\overline{\fR_{l,k}} + \fR'_{l,k})
		= 4\pi\ii \conj^{l+1} F_{l,k}.
	\]
	As for $k = 1$ and $l < L$, we compute at $\tau=0$
% 	and use induction to show that 
% 	\begin{equation}
% 		\frac{d\fF'_{l,1}}{d(\tau^2)} = 
% 		4\pi\ii \conj^{l+1} \Big(
% 			F_{l,1} + \sum_{j=1}^{l-1}\sum_{k=1}^{n_j} 		F_{j,k}
% 		\Big).
% 		\label{eq:dF_{l,1}}
% 	\end{equation}
% Whe $l=1$, 
% \[
% \begin{split}
% \frac{d\fF'_{l,1}}{d(\tau^2)}&=\sum_{k=2}^{n_1}\frac{d\fF'_{l,k}}{d(\tau^2)}+4\pi i\sum_{k=1}^{n_1}\conj^2F_{1,k}\\
% &=-\sum_{k=2}^{n_1}4\pi i\conj^2F_{1,k}+4\pi i\sum_{k=1}^{n_1}\conj^2F_{1,k}\\
% &=4\pi i\conj^2F_{1,1},
% \end{split}
% \]
% so \eqref{eq:dF_{l,1}} holds when $l=1$.  Suppose \eqref{eq:dF_{l,1}} for some $l$ with $1<l<L$.  Then
% \[
% \begin{split}
% \frac{d\fF'_{l+1,1}}{d(\tau^2)}=&-\sum_{k=2}^{n_{l+1}}\frac{d\fF'_{l+1,k}}{d(\tau^2)}-\sum_{k=1}^{n_l}\conj\left(\frac{d\fF'_{l,k}}{d(\tau^2)}\right)+4\pi i\sum_{k=1}^{n_{l+1}}\conj^{l+2}F_{l+1,k}\\
% =&-\sum_{k=2}^{n_{l+1}}4\pi i\conj^{l+2}F_{l+1,k}-\conj\left(4\pi i\conj^{l+1}\left(F_{l,1}+\sum_{j=1}^{l-1}\sum_{k=1}^{n_j}F_{j,k}\right)\right)\\
% &-\sum_{k=2}^{n_l}\conj\left(4\pi i\conj^{l+1}F_{l,k}\right)+4\pi i\sum_{n_{l+1}}\conj^{l+2}F_{l+1,k}\\
% =&4\pi i\conj^{l+2}\left(F_{l+1,1}+\sum_{j=1}^{l}\sum_{k=1}^{n_j}F_{j,k}\right).
% \end{split}
% \]
% By induction, \eqref{eq:dF_{l,1}} holds for $l=1,2,\ldots,L-1$.
	\begin{align*}
		&\sum_{k=1}^{n_l} \frac{d\fF'_{l,k}}{d(\tau^2)}
		+ \sum_{k=1}^{n_{l-1}} \conj\Big(\frac{d\fF'_{l-1,k}}{d(\tau^2)}\Big)\\
		=& -4\pi\ii\frac{\partial}{\partial\tau^2}\bigg(\sum_{k=1}^{n_l} \alpha_{l,k} -
		\sum_{k=1}^{n_{l-1}} \alpha_{l-1,k}\bigg)\\
		&+4\pi\sigma_{l,0}\frac{\partial}{\partial\tau^2}\bigg(\sum_{k=1}^{n_l} \gamma_{l,k} - \sum_{k=1}^{n_{l-1}} \gamma_{l-1,k}\bigg)
		&& \text{(because $\sigma_{l-1,0} = - \sigma_{l,0}$)}\\
		&+2\pi\ii\bigg(\sum_{k=1}^{n_l} \fR'_{l,k}-\sum_{k=1}^{n_{l-1}} \overline{\fR'_{l-1,k}}\bigg)\\
		=& 4\pi\ii\frac{\partial}{\partial\tau^2}(\alpha_{l,0} + \alpha_{l,\infty})
		-4\pi\sigma_{l,0}\frac{\partial}{\partial\tau^2}(\gamma_{l,0}+\gamma_{l,\infty})
		&& \text{(by~\eqref{eq:resalpha} and~\eqref{eq:resgamma})}\\
		&+2\pi\ii\bigg(
			\sum_{k=1}^{n_l} \fR'_{l,k}
			+ \sum_{k=1}^{n_l} \overline{\fR_{l,k}} + \beta_{l,0}^2 - \beta_{l,\infty}^2
	  \bigg)
		&& \text{(by~\eqref{eq:residuephi2})}\\
		=& 2\pi\ii\Big(\beta_{l,\infty}^2 - \beta_{l,0}^2
			+\sum_{k=1}^{n_l}(\overline{\fR_{l,k}} + \fR'_{l,k})
		+\beta_{l,0}^2 - \beta_{l,\infty}^2\Big)
		&& \text{(by~\eqref{eq:fixalphagamma})}\\
		=& 4\pi\ii \sum_{k=1}^{n_l} \conj^{l+1} F_{l,k}.
	\end{align*}
	Then
	\begin{align*}
		\sum_{k=1}^{n_l} \frac{d\fF'_{l,k}}{d(\tau^2)}
		=&\frac{d\fF'_{l,1}}{d(\tau^2)} + 4\pi\ii \conj^{l+1} \sum_{k=2}^{n_l} F_{l,k}\\
		=&(-\conj)^l \sum_{m=1}^l (-\conj)^m \bigg(\sum_{k=1}^{n_m} \frac{d\fF'_{m,k}}{d(\tau^2)}
		+ \sum_{k=1}^{n_{m-1}} \conj\Big(\frac{d\fF'_{m-1,k}}{d(\tau^2)}\Big)\bigg)\\
		=&(-\conj)^l \sum_{m=1}^{l}(-\conj)^m \Big(4\pi\ii\sum_{k=1}^{n_m} \conj^{m+1} F_{m,k}\Big)\\
		=&4\pi\ii\conj^{l+1}\sum_{m=1}^{l}\sum_{k=1}^{n_m} F_{m,k},
	\end{align*}
	so $\widetilde\fF'_{l,1} := \tau^{-2} \fF'_{l,1}$ extends to
	$\tau=0$ with the value
	\[
		\frac{d\fF'_{l,1}}{d(\tau^2)} = 
		4\pi\ii \conj^{l+1} \Big(
			F_{l,1} + \sum_{m=1}^{l-1}\sum_{k=1}^{n_m} F_{m,k}
		\Big).
	\]
\end{proof}

Therefore, if $(q,\dot\theta)$ is balanced, $\widetilde\fF'=0$ is solved at
$\tau=0$.  Recall that we normalize the complex scaling on $\C^\times_1$ by
fixing $p_{1,1}$.  If $(q,\dot\theta)$ is rigid, because $\Theta_2 = \sum
F_{l,k} = 0$ independent of $p$, the partial derivative of
$(\widetilde\fF')_{(l,k) \ne (L,1)}$ with respect to $(p_{l,k})_{(l,k) \ne
(1,1)}$ is an isomorphism from $\C^{N-1}$ to $\C^{N-1}$.  The following
proposition then follows by the Implicit Function Theorem.

\begin{proposition}\label{prop:p}
	Assume that the parameters $t_{l,k}$, $\alpha_{l,k}$, $\beta_{l,k}$,
	$\gamma_{l,k}$, $p'_{l,k}$ are given by Propositions~\ref{prop:talphagamma},
	\ref{prop:beta}, and \ref{prop:pprime}.  Assume further that the central
	values $q_{l,k} = \conj^l\cent p_{l,k}$ and $\dot\theta_h = \cent\beta_h$
	form a balanced and rigid configuration $(q, \dot\theta)$.  Then for
	$(\tau,\beta)$ in a neighborhood of $(0,\dot\theta)$ that
	solves~\eqref{eq:endbalanceb} and~\eqref{eq:endbalancea}, there exists values
	for $p_{l,k}$, unique up to a complex scaling, depending smoothly on $\tau$
	and $(\beta_h)_{h \in \gH}$, such that $p_{l,k}(0,\dot\theta) = \cent
	p_{l,k}$ and the conformality condition~\eqref{eq:conform3} is solved.
\end{proposition}

\subsection{Embeddedness}

It remains to prove that
\begin{proposition}
	The minimal immersion given by the Weierstrass parameterization is regular
	and embedded.
\end{proposition}

\begin{proof}
	The immersion is regular if $|\Phi_1|^2 + |\Phi_2|^2 + |\Phi_3|^2 > 0$.  This
	is easily verified on $U_\delta$.  On the necks and the ends, the
	regularity follows if we prove that $\widetilde\Phi_2$ has no zeros outside
	$U_\delta$.  At $\tau=0$, $\widetilde\Phi_2$ has $n_l + n_{l-1}+2$ poles on
	$\hat\C_l$, hence $n_l + n_{l-1}$ zeros.  By taking $\delta$ sufficiently
	small, we may assume that all these zeros lie in $U_{l,\delta}$.  By
	continuity, $\widetilde\Phi_2$ has $n_l + n_{l-2}$ zeros in $U_{l,\delta}$
	also for $\tau$ sufficiently small.  But for $\tau \ne 0$, $\widetilde\Phi_2$
	is meromorphic on a Riemann surface $\Sigma_\tau$ of genus $g = N-L$ and has
	$2L+2$ simple poles, hence has $2(N-L) - 2 + 2L + 2 = 2N$ zeros.  So
	$\widetilde\Phi_2$ has no further zeros, in $\Sigma_t$, in particular not
	outside $U_\delta$.

	We now prove that the immersion
	\[
		z\mapsto	\re \int^z (\Phi_1, \widetilde\Phi_2, \Phi_3)
	\]
	is an embedding, and the limit positions of the necks are as prescribed.

	On $U_{l,\delta}$, the Gauss map $G = -(\Phi_1 + \ii \Phi_2)/\Phi_3$ converges
	to $\ii\sigma_{l,0}$, so the immersion is locally a graph over the $xz$-plane.
	Fix an orientation $\sigma_{1,0} = -1$, then up to translations, we have
	\[
		\lim_{\tau \to 0} \Big(\re \int^z \Phi_1 + \ii \re \int^z \Phi_3\Big)
		= \conj^l(\ln z) + 2m\pi\ii
	\]
	where $m$ depends on the integral path, and
	\[
		\lim_{\tau \to 0} \re \int^z \widetilde\Phi_2 = \re \int^z (\conj^*)^l \psi_l
		=: \Psi_l(\conj^l z)
	\]
	which is well defined for $z \in U_{l,\delta}$ because the residues of $\psi_l$
	are all real.

	With a change of variable $z \mapsto \ln z$, we see that the immersion
	restricted to $U_{l,\delta}$ converges to a periodic graph over the
	$xz$-planes, defined within bounded $x$-coordinate and away from the points
	$\ln q_{l,k} + 2m\pi\ii$, and the period is $2\pi\ii$.  Here, again, we
	identified the $xz$-plane with the complex plane.

	This graph must be included in a slab parallel to the $xz$-plane with bounded
	thickness.  We have seen from the integration along $B_k$ that the distance
	between adjacent slabs is of the order $\mathcal{O}(\ln \tau)$.  So the slabs
	are disjoint for $\tau$ sufficiently small.

	As for the necks and ends, note that there exists $Y > 0$ such that
	$\Psi_l^{-1}([-Y,Y])$ is bounded by $n_l+n_{l-1}+2$ convex curves.  After the
	change of variable $z \mapsto \ln z$, all but two of these curves remain
	convex; those around $0$ and $\infty$ become periodic infinite curves.  If
	$Y$ is chosen sufficiently large, there exists $X > 0$ independent of $l$
	such that the curves $|z| = \exp(\pm X)$ are included in
	$\Psi_l^{-1}([-Y,Y])$ for every $1 \le l \le L+1$.  After the change of
	variable $z \mapsto \ln z$, these curves become curves with $\re z = \pm X$.

	Hence for $\tau$ sufficiently small, we may find $Y^+_l$ and $Y^-_l$, with
	$Y^-_l < Y^+_l < Y^-_{l+1}$, and $X>0$, such that
	\begin{itemize}
		\item The immersion with $Y^-_l < y < Y^+_l$ and $-X < x < X$ is a graph
			bounded by $n_l + n_{l-1}$ planar convex curves parallel to the $xz$-plane,
			and two periodic planar infinite curves parallel to the $yz$-plane.

		\item The immersion with $Y^+_l < y < Y^-_{l+1}$ and $-X < x < X$
			consists of annuli, each bounded by two planar convex curves parallel to
			the $xz$-plane.  These annuli are disjoint and, by a Theorem of
			Schiffman~\cite{schiffman1956}, all embedded.

		\item The immersion with $|x| > X$ are ends, i.e.\ graphs over
			vertical half-planes, extending in the direction $(-1, -\dot\theta_{l,0})$ and
			$(+1, -\dot\theta_{l,\infty})$, $1 \le l \le L+1$.  If the
			inequality~\eqref{eq:embedcond} is satisfied, these graphs are disjoint.
	\end{itemize}
	This finishes the proof of embeddedness.
\end{proof}

\bibliography{References}

\begin{thebibliography}{dSRB10}

\bibitem[CF22]{chen2022}
Hao Chen and Daniel Freese.
\newblock Helicoids and vortices.
\newblock {\em Proc. R. Soc. A}, 478(2267):20220431, 2022.

\bibitem[Che21]{saddle2}
Hao Chen.
\newblock {Gluing Karcher-Scherk saddle towers II: Singly periodic minimal
  surfaces}, 2021.
\newblock preprint, \texttt{arXiv:2107.06957}.

\bibitem[Con17a]{connor2017}
P.~Connor.
\newblock A note on balance equations for doubly periodic minimal surfaces.
\newblock {\em Math. J. Okayama Univ.}, 59:117--130, 2017.

\bibitem[Con17b]{connor2017b}
P.~Connor.
\newblock A note on special polynomials and minimal surfaces.
\newblock {\em Houston J. Math.}, 43:79--88, 2017.

\bibitem[CT21]{saddle1}
Hao Chen and Martin Traizet.
\newblock {Gluing Karcher-Scherk saddle towers I: Triply periodic minimal
  surfaces}, 2021.
\newblock preprint, \texttt{arXiv:2103.15676}.

\bibitem[CW12]{connor2012}
P.~Connor and M.~Weber.
\newblock The construction of doubly periodic minimal surfaces via balance
  equations.
\newblock {\em Amer. J. Math.}, 134:1275--1301, 2012.

\bibitem[DJJ13]{dominici2013}
D.~Dominici, S.~J. Johnston, and K.~Jordaan.
\newblock Real zeros of {${}_2F_1$} hypergeometric polynomials.
\newblock {\em J. Comput. Appl. Math.}, 247:152--161, 2013.

\bibitem[dSRB10]{dasilva2010}
M.~F. da~Silva and V.~Ramos~Batista.
\newblock Scherk saddle towers of genus two in {$\mathbb R^3$}.
\newblock {\em Geom. Dedicata}, 149:59--71, 2010.

\bibitem[Hei78]{heine1878}
Heinrich~Eduard Heine.
\newblock {\em Handbuch der Kugelfunctionen, Theorie und Anwendungen, Bd. I}.
\newblock G. Reimer, Berlin, second edition, 1878.

\bibitem[HLB14]{hancco2014}
A.~J.~Yucra Hancco, G.~A. Lobos, and V.~Ramos Batista.
\newblock Explicit minimal scherk saddle towers of arbitrary even genera in
  $\mathbb{R}^3$.
\newblock {\em Publ. Mat.}, 58(2):445--468, 2014.

\bibitem[HMR09]{hauswirth2009}
L.~Hauswirth, F.~Morabito, and M.~Rodr\'{i}guez.
\newblock An end-to-end construction for singly periodic minimal surfaces.
\newblock {\em Pacific Journal of Math}, 241(1):1--61, 2009.

\bibitem[Kar88]{karcher1988}
H.~Karcher.
\newblock Embedded minimal surfaces derived from {S}cherk's examples.
\newblock {\em Manuscripta Math.}, 62(1):83--114, 1988.

\bibitem[Li12]{li2012}
Kevin Li.
\newblock {\em Singly-periodic minimal surfaces with Scherk ends near parallel
  planes}.
\newblock PhD thesis, Indiana University, Bloomington, 2012.

\bibitem[Mar66]{marden1966}
Morris Marden.
\newblock {\em Geometry of polynomials}.
\newblock Mathematical Surveys, No. 3. American Mathematical Society,
  Providence, R.I., second edition, 1966.

\bibitem[MRB06]{martin2006}
Francisco Mart\'{\i}n and Val\'{e}rio Ramos~Batista.
\newblock The embedded singly periodic {S}cherk-{C}osta surfaces.
\newblock {\em Math. Ann.}, 336(1):155--189, 2006.

\bibitem[MW07]{meeks2007}
William~H. Meeks, III and Michael Wolf.
\newblock Minimal surfaces with the area growth of two planes: the case of
  infinite symmetry.
\newblock {\em J. Amer. Math. Soc.}, 20(2):441--465, 2007.

\bibitem[Sch35]{scherk1835}
H.~F. Scherk.
\newblock Bemerkungen \"{u}ber die kleinste {F}l\"{a}che innerhalb gegebener
  {G}renzen.
\newblock {\em J. Reine Angew. Math.}, 13:185--208, 1835.

\bibitem[Shi56]{schiffman1956}
Max Shiffman.
\newblock On surfaces of stationary area bounded by two circles, or convex
  curves, in parallel planes.
\newblock {\em Ann. of Math. (2)}, 63:77--90, 1956.

\bibitem[Tra96]{traizet1996}
Martin Traizet.
\newblock Construction de surfaces minimales en recollant des surfaces de
  {S}cherk.
\newblock {\em Ann. Inst. Fourier (Grenoble)}, 46(5):1385--1442, 1996.

\bibitem[Tra01]{traizet2001}
Martin Traizet.
\newblock Weierstrass representation of some simply-periodic minimal surfaces.
\newblock {\em Ann. Global Anal. Geom.}, 20(1):77--101, 2001.

\bibitem[Tra02a]{traizet2002b}
Martin Traizet.
\newblock Adding handles to {R}iemann's minimal surfaces.
\newblock {\em J. Inst. Math. Jussieu}, 1(1):145--174, 2002.

\bibitem[Tra02b]{traizet2002}
Martin Traizet.
\newblock An embedded minimal surface with no symmetries.
\newblock {\em J. Diff. Geom.}, 60:103--153, 2002.

\bibitem[Tra08]{traizet2008}
Martin Traizet.
\newblock On the genus of triply periodic minimal surfaces.
\newblock {\em J. Diff. Geom.}, 79(2):243--275, 2008.

\bibitem[Tra13]{traizet2013}
Martin Traizet.
\newblock Opening infinitely many nodes.
\newblock {\em J. Reine Angew. Math.}, 684:165--186, 2013.

\bibitem[TW05]{traizet2005}
M.~Traizet and M.~Weber.
\newblock Hermite polynomials and helicoidal minimal surfaces.
\newblock {\em Inv. Math.}, 161:113--149, 2005.

\end{thebibliography}
\bibliographystyle{alpha}
\end{document}